\documentclass[12pt]{article}
 \usepackage{amsmath,amsthm,amssymb}
\usepackage{graphicx}
\usepackage{wrapfig}
\usepackage{pifont}
\usepackage{color}

\begin{document}

\theoremstyle{plain}
\newtheorem{Thm}{Theorem}[section]
\newtheorem{Prop}[Thm]{Proposition}
\newtheorem{Pp}[Thm]{P}
\newtheorem{Defi}[Thm]{Definition}
\newtheorem{Cor}[Thm]{Corollary}
\newtheorem{Lem}[Thm]{Lemma}

\makeatletter
\renewcommand{\theequation}{%
\thesection.\arabic{equation}}
\@addtoreset{equation}{section}
\makeatother

\setlength{\baselineskip}{16pt}

\def\textmc{\rm}
\def\({(\!(}
\def\){)\!)}
\def\R{\mathbb{R}}
\def\Z{\mathbb{Z}}
\def\N{\mathbb{N}}
\def\H{\mathbb{H}}
\def\C{\mathbb{C}}
\def\D{\mathbb{D}}
\def\Q{\mathbb{Q}}
\def\T{\mathbb{T}}
\def\E{{\bf E}}
\def\pr{{\bf P}}
\def\M{{\cal M}}     
\def\F{{\cal F}}
\def\G{{\cal G}}
\def\cD{{\cal D}}
\def\X{{\cal X}}
\def\A{{\cal A}}
\def\B{{\cal B}}
\def\L{{\cal L}}
\def\a{\alpha}
\def\b{\beta}
\def\e{\varepsilon}
\def\de{\delta}
\def\ga{\gamma}
\def\k{\kappa}
\def\la{\lambda}
\def\fa{\varphi}
\def\th{\theta}
\def\si{\sigma}
\def\t{\tau}
\def\om{\omega}
\def\De{\Delta}
\def\Ga{\Gamma}
\def\La{\Lambda}
\def\Om{\Omega}
\def\Th{\Theta}
\def\lan{\langle}
\def\ran{\rangle}
\def\lbr{\left(}
\def\rbr{\right)}
\def\v2{\vskip2mm}
\def\n{\noindent}
\def\pf{{\it Proof.~}}
\def\1{{\bf 1}}

\def\Tr{\operatorname{Tr}}
\def\quadd{\qquad\qquad}
\def\n{\noindent}
\def\beq{\begin{eqnarray*}}
\def\eeq{\end{eqnarray*}}
\def\supp{\mbox{supp}}
\def\beqn{\begin{equation}}
\def\eeqn{\end{equation}}

%%%%

\begin{center}
        \end{center}
                \begin{center}
  {\Large  Boundary behaviour  of  RW's on planar graphs\\
      and convergence of LERW to  chordal SLE$_2$ }
\\
\vskip4mm
{K\^ohei UCHIYAMA} \\
\vskip2mm
{Department of Mathematics, Tokyo Institute of Technology} \\
{Oh-okayama, Meguro Tokyo 152-8551\\
e-mail: \,uchiyama@math.titech.ac.jp}
\end{center}

\vskip8mm
\n
{\it running head}:    RW's on planar graphs\\
\n
{\it key words}: random walks  on  a planar graph;   loop-erased random walk;  
stochastic Loewner evolution; chordal  SLE

\n

\vskip6mm
\begin{abstract}   This paper  concerns  a random walk  on  a planar graph  and presents  certain estimates concerning the harmonic measures for the walk in a grid domain   which estimates are useful for showing  the convergence of  a LERW (loop-erased random walk) to an  SLE (stochastic Loewner evolution). 
We assume that  the walk started at a  fixed vertex of the graph satisfies the invariance principle as in
 Yadin and Yehudayoff \cite{YY} in which the convergence of LERW to a radial SLE is established  in this setting.  
  Our main concern is chordal case, where a random walk is started at a boundary vertex  of a  simply connected grid domain and conditioned to exit it through another boundary vertex specified in advance.  The primary contribution of the present  paper is an   
estimate, which states   that the excursion of the conditioned  walk  leaves an intrinsic neighborhood of its initial point not \lq along' the boundary but through an intrinsic interior of the domain with high probability. Based on this result we give a proof for the convergence to the chordal  SLE, a result that has recently been 
proved by  Suzuki \cite{Su} under  an analyticity assumption on  the  boundary of the domain arising in the limit.
\end{abstract}

%%%%%%%%%%%
%%%%%%%%%
\newpage
\begin{center} Contents
        \end{center}
        
1. \;Introduction.   

2. \;Planar graph and Hypothesis (H).

\quad 2.1. Planar graph and random walks on it.

\quad  2.2.  The metric of the path space and hypothesis (H).

\quad 2.3. A consequence of Hypothesis (H).

\quad 2.4. Supplement for the proof of Proposition \ref{YY4.1}.
 
3. \;Preliminary lemmas of geometric nature.     
 
 \quad 3.1. Elementary bounds on distortion under $\fa_D$. 
   
   \quad 3.2.   Domains   $ B_r(x), U_\de, Q(x)$.  \; 
   
   \quad  3.3.   Construction of a domain $U$.

4. \;Estimates of hitting distributions of random walks.

\quad 4.1.  Simple properties  of the planar graph. 

  \quad  4.2. Starting near the boundary (unconditional case).

 \quad   4.3.  Starting near the boundary (conditional case). 
 
  \quad 4.4. Hitting distribution of $\partial D$.  
  
  \quad 4.5. Poisson kernel approximation.

5. \;Convergence of LERW to chordal SLE$_2$.

\quad 5.1. Chordal Loewner chains for a simple curve in $\H$. \; 

\quad  5.2. Chordal Loewner chains in a simply connected domain. 

\quad 5.3. Convergence of driving function. 

\quad 5.4. Uniform convergence.

\vskip4mm
\section{Introduction}

This paper  concerns  a random walk  on  a planar graph  imbedded in the plane and provides  certain estimates concerning the harmonic measures for the walk in a domain. The estimates obtained  are used for showing  the convergence of  a loop-erased random walk (LERW) to a stochastic  Loewner evolution (SLE).  Our essential hypothesis is  that the walk started at one fixed vertex of the planar graph  satisfies  the 
invariance principle (as in  \cite{YY}): properly scaled trajectory of it weakly converges to that of the planar Brownian motion with respect to a metric which disregards the difference of time parametrization. We do not assume the symmetry of the random walk, while  the planarity of the graph plays an essential role  
 as in \cite{YY}.

The loop-erased  random walk is a process obtained by  erasing loops one by one from a random walk on a graph in  chronological order. It   was  introduced by Lawler \cite{L_D} as a version of self-avoiding random walk focusing the central limit behavior (functional limit theorem) in dimensions $\geq 4$,  and  there have appeared many works studying  various aspects of it (\cite{D}, \cite{GB}, \cite{L2}, \cite{K1},\cite{K2}, \cite{M}, etc.).  %It turns out that  in two dimensions  its scaling limit of LERW is  SLE$_2$. 

The stochastic---or Schramm---Loewner evolutions (SLE$_\k)$  are a family of  random trajectories  obtained as a solution of the Loewner differential equation (in the complex plane) driven by the process $\sqrt{\k}\, W(t)$, where  $W$ is a 
one-dimensional standard Brownian motion and $\k$ is a positive parameter.
The  SLE's  are introduced by Schramm in \cite{S},  in which it is  conjectured,  based on   ample  evidence provided by foregoing works and  his paper itself,  that the scaling limit of LERW on $\Z^2$ must be SLE$_2$. 

   This   conjecture by Schramm is proved by Lawler, Schramm and Werner  \cite{LSW} where  a scaled LERW
on some regular  lattices is shown    to converge to a {\it radial}  SLE$_2$, a version of SLE path that starts from an interior point of the domain and ends up at a boundary point of it. 
Dapeng Zhan \cite{Z} have studied  LERW's on the square lattice   in a multiply connected domain 
and proved the existence of their scaling limit;
in the case of a simply connected domain in particular, he  has proved the convergence to a {\it chordal} SLE$_2$, another version of SLE path that travels from  boundary  to boundary.
Yadin and Yehudayoff \cite{YY} extend the result of   \cite{LSW},  the convergence of LERW to a radial SLE to that
for the natural random walks on planar graphs under  a  natural setting (the same as ours) of the problem. 
Recently Suzuki \cite{Su}  have obtained a chordal version of their result: the LERW  in a  simply connected domain conditioned to connect two boundary vertices  converges to a chordal SLE$_2$ curve  in a  setting similar to  \cite{YY} under the assumptions (1)  the invariance principle holds uniformly for starting points of the walk and more seriously (2) the boundary of the domain is locally analytic at the  starting  boundary point   of the random walk from whose trajectory (or rather its time reversal)   the LERW is derived.

In this paper we are  concerned with the chordal case of LERW  on a planar graph as  in \cite{Su}.   For the chordal case of LERW  we need unlike the radial  case to deal with an excursion of random walk (a random walk path  in a domain connecting a boundary point  with another  one),
 and   we  are forced to estimate the harmonic measure of the random walk started at a vertex on (or near) the boundary and conditioned to exit the domain through another boundary vertex that is specified in advance.   The SLE$_2$  curve  is  conformally invariant, to which we intend to  show the LERW's obtained from such excursions scaled by the sizes of the domains converge.
The approximation   must accordingly be effected  uniformly  for the domain  as far as its  sizes (measured by the  inner radius with respect to an appropriately chosen point) is large enough and in order to obtain such uniformity  we wish to find  a certain estimate concerning  the distribution of \lq entrance' of the excursion  into the substantial interior of the domain. The primary contribution of the present  paper is such an   
estimate  (Proposition \ref{prop2})  that  states   that the excursion leaves an intrinsic neighborhood of its initial point not \lq along' the boundary but through an intrinsic interior of the domain  with high probability.
The result plays a key role  in our verification of   the convergence of LERW to  the chordal  SLE:  it refines the  result of  \cite{Su} by removing the second extra assumption mentioned above.

%Our basic assumption ((H) in Section 3)  is that the  random walk started at a  prescribed vertex  behaves like Brownian motion,  as already  mentioned.     Under  this  assumption we show that the  invariance principle then  actually holds uniformly w.r.t.  the starting position,  which  makes the analysis concerning the random walk easier.  Actually it simplifies  not a little of 
%the arguments made by Yadin and Yehudayoff \cite{YY} of which  the essential result is to establish an approximation formula for the ratio of  the discrete Poisson kernels  of a domain  by the corresponding continuous  version.   This result of them,  of which we  include a  proof, is also  applied to the present issue on  the convergence of LERW to  the chordal  SLE.

\v2
\v2

%%%%%%%%%  
The rest of the paper is organized as follows.  %In Section 2 we introduce a simple coupling argument of two independent random walks  with distinct starting points.  This  is the  main tool for  the proof of the uniformity of invariance principle, which we give
 In Section 2   we introduce our planar graph and the  random walk on it,  and then state  Hypothesis (H), our basic assumption in this paper,  and derive a  lemma from it by using a key result of \cite{YY}.     In Section 3  we discuss some geometric properties of a conformal map from a simply connected domain onto the unit disc  and derive  preliminary facts used in the next section.   Section 4 consists of 5 subsections and provides  various 
 estimates concerning  the hitting distribution  of the random walk; Proposition \ref{prop2}, the primary result of the paper, is    proved  in the third subsection of it.
   %  Poisson approximation is included in the last subsection.
 Section 5 concerns the convergence of  a LERW on the planar graph  to a chordal SLE$_2$. We break this section into four  subsections.  
A brief review  of the chordal SLE in $\H$  is provided in Section 5.1. In Section 5.2  we present an account of  the chordal SLE in a simply connected domain together with some facts concerning it. The statement of  the  convergence result  together with its abridged proof is given in the last two subsections.

%%%%%\section{Uniform invariance principle} 
\section{Planar graph and Hypothesis (H)}

Here we introduce our random walk on a planar graph as well as Hypothesis (H), and discuss on the fundamental facts on  the harmonic measures of the scaled walk. The setting is the same as in  \cite{YY}.    We formulate a key idea used in it as  a lemma that is convenient  to apply in the rest of the present paper.

\v2\n
{\sc 2.1. \;Planar graph and random walks on it.}

\v2
For any $x, y \in \C $, we write $[x, y]$ for $\{(1-t)x+ty :0\leq t\leq 1\} $,  the  directed  line segment emanating from   $x$ and ending at $y$. (The same square brackets is used to designate a closed interval of $\R$, but  this will cause  no confusion.)   
Let $V$ be a countable subset of $\C$. Let $p:V\times V\mapsto [0,1]$ be such that
$$\sum_{v\in V} p(u,v)=1$$
and put $E=\{ (u, v) : p(u, v)>0 \} $. A pair  $(u,v)\in E$ is identified with the directed segment $[u,v]$ and  is called an  edge.   The pair $G=(V, E)$  may be considered to be a graph of directed edges $[u,v]\in E$ with weight $p(u,v)$.  
%We assume that $\sum _ {v \in V}E(u, v) < \infty$ for every $u \in V$, and put
%\[ p(u,v)=\frac{E(u, v)}{\sum _{w \in V} E(u,w)} .\]
We are concerned with the  Markov chain  whose transition probability  is  $p(u,v)$.

In this paper, we assume that the graph $G$ satisfies the following properties.
\begin{enumerate}
\item
 $G$ is a planar graph, namely  any two  edges  are disjoint unless they have at least one  common  endpoint.
\item
For any compact set $K \subset \C$,  $\sharp (K\cap V)<\infty$.  ($\sharp$ designates the cardinality of a set.)  

\item
The Markov chain $(S_k)_{k=0}^\infty$  on $V$ with transition probability $p(u,v)$ is irreducible.
%\item
%There exists an invariant measure $\pi$ on $V$ for $S_n$ such that
%\[0<\pi (u)<\infty ,\quad u \in V.\]

\end{enumerate}
For simplicity we further suppose that  $0\in V$.

We suppose that for each $v\in V$ there is given a random sequence $(S^v_k)_{k=0}^\infty \subset V$ defined on  a probability space $(\Om, \F,P)$ that constitutes a Markov chain  with transition probability $p(u,v)$   such that  $S^v_0=v$ and that any two sequences with distinct initial vertices are independent. We denote the linear interpolation of $(S^v_k)$ by $(\bar S^{\,v}_t)_{t\geq 0}$: $\bar S^{\,v}_t$ travels along the edge $[S_n^v, S_{n+1}^v]$  with unit speed for $t\in [n,n+1]$. For any  subset $U$ of $\C$ and for $v\in V$  
we define 
the first exit time of $ S^v$ from $U$, denoted by $\tau_U$, as  the least positive integer  $k$ such that the segment $[S^v_{k-1}, S^v_{k}]$ contains a point  of $\C\setminus U$: 
%$$\bar\tau_U = \inf \{t\geq 1: \bar S^v_t \notin U\}$$ 
%and put $\tau_U=  \lceil  \bar\tau_U \rceil $, where $\lceil a \rceil $ denotes the least integer that is not smaller than $a$. In other words $\tau_U$ is
$$ \tau _U= \inf \{ k\geq 1 : [S_{k-1}^v, S_{k}^v ] \setminus U\neq \emptyset\}. 
$$
If a set  $K$ intersects the segment $[S^v_{k-1}, S^v_k]$ for $k=\tau_U$, then  $S^v$ is said to  {\it exit $U$ through}   $K$.  ($K$ will be contained in the complement of $U$ in our use.)
 %(For  the case  $v\in \partial U$ see Section 4, where  we shall be interested in the walk started at a vertex outside $U$.) 
  (In the case $v\in \partial U$ we shall modify the definition of $\tau_U$ in Section 4.3, until that time we shall not encounter the situation where the modification is needed.)
Obviously  $\tau_U$ depends on $v$ and we sometimes indicate this dependence by writing $\tau^v_{U}$ but usually do not when   it is  clear from the context. We also suppose that the standard Brownian 
motion on $\C$ is defined on  $(\Om, \F,P)$ and denote it by $W_t$ and write $W^z_t =z + W_t$, $z\in \C$ (so that $W_0^z=z$). The first exit time of $W^z_t$ from $U$ will be denoted by $\tau^W_U$: $\tau^W_U=\inf \{t>0: W^z_t \notin U\}$.
% (we shall explicitly state which process $\tau_U$ is associated with   in the situation where ambiguity arises).
 
For a set $A\subset \C$  we write $V(A)$ for $V\cap A$ and denote by $G(A)$ the subgraph 
of $G$ of which the vertex set is $V(A)$ and the edges  are those $(u,v)\in E$ such that $[u,v]\subset A$. 
A path in  $G( A)$ is a finite sequence $u_0,\ldots, u_n$ such that $(u_{k-1}, u_{k})$ is an edge of $G(A)$ for each $k=1,\ldots, n$. If $v\in V(U)$, $\tau_U$ may agrees with the first time  when  $S_n^v$, considered to be a walk on $G$,  exits the subgraph $G(U)$.

%%%%%%%
%%%%
%If a linear interpolation of $S^\rho (n) $ started at $0$ and stopped on exiting ${D}$
%converges weakly to a two-dimensional standard Brownian motion started at $0$ and stopped on exiting ${ D}$
%as $\rho \rightarrow \infty$ with respect to $d_\mathcal{U}$, then we say
%$S_\rho (n) $ satisfies invariance principle relative to $D$.

%2.2
\v2\n
{\sc 2.2. \; The metric of the path space and hypothesis (H).}
\v2

Denote by   ${\cal C}^*$ the  space of finite continuous plane curves. Here it is understood that a curve  (also called path)  is  oriented and   represented by  $\ga \in C[0,\tau]$, a continuous map of a finite interval $[0,\tau]$ into $\C$, but two such  maps are identified if they are transformed to each other by some changes of parametrization that preserve  orientation. 
%a trajectories on $D$ ending at a point of $\partial D$: $\ga\in {\cal C}_*(D)$ is a map of a finite interval $[0,\tau]$ into $\bar D:= D\cup \partial D$ such that $\ga[0,\tau)\subset D$ and $\ga(\tau) \in \partial D$. On ${\cal C}_*(D)\times {\cal C}_*(D)$  define a quasi-metric  $d_D$ 
We consider ${\cal C}^*$ as a metric space with the metric $d_{{\cal U}}^*$ defined as follows: for  $\ga_j \in C[0,\tau_j]$ $(j=1,2)$,
$$d_{\cal U}^* (\ga_1,\ga_2) = \inf_{\chi}\sup_{0\leq t\leq \tau_1}|  \ga_1(t) -    \ga_2 \circ\chi(t)|,$$
where  the infimum is taken over all  homeomorphisms $\chi : [0,\tau_1] \mapsto [0,\tau_2] $ with $\chi(0)=0$. The metric space ${\cal C}^*$  is separable and complete.

%the set of equivalence class ${\cal C}_*(D)/\sim$ is  a separable metric space with the metric  $\rho_D$  function  induces the equivalence relation $\sim$ in ${\cal C}_*(D)$
%Let $W^x_t$ be a standard  planar Brownian motion stated at $x\in \C$.  For a measurable function $F$ of paths in ${\cal C}^*(\D)$, put
%$$\hat F_D(x)= E[ F_D(W^u)]$$
%where $\tau_D =\inf \{t>0: W_t \notin D\}$. Here $F_D(\ga) = F( \ga_{t\wedge  \tau_D})_{t\geq 0})).$
Under this metric we shall consider the convergence of probability measures on ${\cal C}^*$ induced by $\bar S^{v}$ that is stopped on exiting a domain $D$ (i.e., at the time $\bar \tau_D$).  Given  a map $\ga\in C[0,\tau]$ with $\ga(0)\in D$,  that represents an element ${\cal C}^*$,  let $\ga^D$ denote the restriction of $\ga$ on $[0,\tau_D\wedge \tau]$ where $\tau_D= \inf\{t\in [0,\tau]: \ga(t)\notin D\}$ and put  
$$d_D(\ga_1, \ga_2) = d_{\cal U}^* (\ga_1^D,\ga_2^D).$$ 

For $h>0$ we scale our random walk  by $h>0$.  The scaling is plainly given by simply multiplying $S_n^v$ by $h$. Thus the scaled walk   $S^{h,v}_n:= h S^v_n$ ($v\in V$) is   the random walk  on $h V= \{h v: v\in V\}$ started at  $h v$; also write  $\bar S^{h,v}$ for the linear interpolation of  $S^{h,v}$.
%Let  some sequence $h_n>0$  decreasing to $0$ be given,
 We denote the open unit disk by $\D$ and
  suppose that the  walk $S^{0}$ %($0<h\leq 1$) stopped on  exiting $\D$ 
   satisfies invariance principle in the following sense: 
\[
{\rm (H)}\;\;\; \left\{\begin{array}{l}  \mbox{\it the law of the scaled walk $\bar S^{\,h,0}$ stopped on exiting  $\D$ weakly  converges to the} \\ 
\mbox{\it law of Brownian motion $W^0$ stopped on exiting $\D$ as $h\downarrow 0$, where the weak} \\
  \mbox{\it convergence is  relative to the metric $d_{\D}$. }
 \end{array}\right.
\]
In \cite{YY}  it is deduced from this hypothesis that for any  $v\in V$ with $hv\in \D,$ the Markov chain $hS^v$ killed on exiting $\D$ behaves like a Brownian motion killed on $\partial\D$ as $h\downarrow 0$ as far as  the hitting 
distributions are  concerned.  We formulate a consequence of this in the next subsection as  Lemma \ref{Nlem1}.  All applications of (H) in this paper will be  via it or its corollaries.

%2.3
\newpage
\v2\n
{\sc 2.3.\; A consequence of Hypothesis (H)}

\v2
Denote by $ \Gamma_r(a) $ the disc of radius  $r>0$  centered at $a\in \C$: 
\beqn\label{disc0}
 \Gamma_r(a) =\{u\in \C: |u-a|<r \}.
 \eeqn
 Let $r>0$ and $U$ be a domain that contains $\Ga_{r'}\setminus \Ga_r$ for some  $r'>r$.  A (continuous) curve  $\ga\in C[0,\tau)$ is said to {\it encompass} the disc $ \Gamma_r(a) $ in  $U$ %the annulus $ \Gamma_R(a)\setminus  \Gamma_r(a) $
  if there is a pair  $s<s' \leq \tau_\ga$ such that 
 \v2
(1) \; $\ga[0,s']\subset U,\;  \ga[s, s'\,] \subset  U \setminus  \Gamma_r(a), \; \ga(s)=\ga(s') \;$  and 
\v2
(2)\;\; the argument    $\arg(\ga(t) -a)$  continuously  varies  from $\arg(\ga(s) - a)$ either to 

\qquad $\arg(\ga(s) -a) +2\pi$ or to $\arg(\ga(s) -a) -2\pi$ as $t$ increases from $s$ to $s'$.  
\v2\n
Here $\ga [s,s']$ designates the restriction of $\ga$ to $[s,s']$.
Because of the condition (1) the term \lq encompass' entails that encompassing is made before exiting $U$.
A random walk encompasses $\Ga_r(a)$ in $U$ if its linear interpolation does. 
 The following result is essentially Proposition 4.1 of \cite{YY} that we adapt and modify to the present need and notation.
 %Prop2.1
 \begin{Prop}\label{YY4.1} \, For any $\e>0$ and $M\geq 1$ there exists $\eta = \eta(\e, M) >0$ such that for any positive number  $r$ one can choose $h_0=h_0(\e,\la, M)>0$ so that if $0< h<h_0$ and  $a\in M\D$, then for all $u\in V(h^{-1}\Ga_{\eta\la}(a))$, 
\[
P\Big[ \bar S^{h,u}\; \mbox{{\rm encompasses}} \; \Gamma_{\eta r}(a) \; \mbox{{\rm in}}\; \Gamma_r(a) \, \Big] >1-\e.
\] 
 \end{Prop}
 \v2\n
 \pf\, In our proposition the statement in \cite[Proposition 4.1]{YY} is modified in two ways. Firstly it is stated for a simply connected domain $D$ by means of   a conformal map $\fa_D$. Our proposition is specialized to the case $D= (M+1)\D$. Secondly the choice of $h_0$ may depend on $a$ in \cite{YY}, while it does not in ours. This independence of $h_0$ from $a$  is verified by
 examining the proof  in \cite{YY}. (We provide  more details in Appendix for the latter.)  \qed
 \v2
 In (H) the domain $\D$ plays no intrinsic role: it may be replaced by any bounded domain containing the origin because  of the scaling property of Brownian motion, hence   the assertion for  $M>1$ follows from that for $M=1$, so we usually  state  results  only  for the case $M=1$   in the sequel.  
%For the same reason  it suffices to obtain (\ref{eq_thm1})  for  $z\in (1-\e)\D$ instead of $z\in \D$.

\v2
We formulate Lemma \ref{Nlem1}, mentioned previously, in terms  not of the scaled walk $S^{h,u}$ but of $S^u$ itself for convenience of later applications. 
Let $U$ be a domain of $\C$ and $K$ a  compact set that is contained in the complement of $ U$. 
For 
$u\in V(U)$  and  $0<r< \rho$,  define
$$q_{K,U}(u) = P[S^{u}  \;\mbox{exits $U$  through} \;K]$$
and% if $\Ga_r(u)\subset U$, 
$$q_{K,U}^{(0)}(u;r) = P[S^{0}\circ\th_{\sigma(\Ga_{r}(u))} \;\mbox{exits $U$  through} \;K\,|\, \sigma_{\Ga_{r}(u))}<\tau_{\rho\D}],$$
where $\sigma_B$ (or $\sigma(B)$) is the first hitting time of  a set $B$ by $S^u$: $\sigma_B = \sigma_{\C\setminus B}$ and   $\th_\sigma$  denotes the usual shift operator so that $ S^{0}\circ \th_\sigma$ is the  the Markov chain  $(S^{0}_{n+\sigma})_{n=0,1, 2,\ldots}$.

\v2
%Lem2.2
\begin{Lem} \label{Nlem1}  Let $\eta =\eta(\e,1)$ be as in Proposition \ref{YY4.1}.  Then for any $\e>0$ and  $\la>1$ one can choose  $R>1$ independently   of $U$ and  $K$ so that  if $\rho \geq R$, then for all $u\in V(\D)$ with 
${\rm dist}(u,\partial U)> r$,
$$(1-\e)q_{K,U}^{(0)}(u; \eta \la)  < q_{K,U}(u)< (1-\e)^{-1}q_{K,U}^{(0)}(u; \eta \la) .$$ 
 \end{Lem}
 \v2\n
 \pf\, We omit $K, U$ from $q_{K,U}^{(0)}(u;r)$ and  $q_{K,U}(u)$.      %Take $\a>0$ small enough that $(1-\a)^{-1}<1+\e$  and    Let $\eta$ be a positive number   corresponding to $\a$ in  Proposition \ref{YY4.1}. 
 Since by strong Markov property $q^{(0)}(u; \eta\la)$ is a convex sum of $q(v)$ over $v\in \Ga_{\eta\la}(u)$ with $P[S^{0}_{\sigma(\Ga_{\eta\la}(u))} = v]>0$, there exist two sites $u^*$ and $u_*$ in $V(\Ga_{\eta\la}(u))$ such that 
 \beqn\label{neq1}
 q(u_*) \leq q^{(0)}(u;\eta\la) \leq q(u^*).
 \eeqn 
  By the
 maximal principle applied to the  stopped chain $(S^{v}_{n\wedge \tau(U)})_{ n=0, 1, \ldots}$ there exists a path $\ga$ in $G(U)$ such that $q(v)\geq q(u^*)$ for $v\in \ga$ and $\ga$ connects $u^*$ with a site outside   $\Ga_{\la}(u)$. Then by Proposition \ref{YY4.1}  the walk  $S^{u}$ intersects $\ga$ before exiting $\Ga_{\eta \la}(u)$ with a probability larger than  $1-\e$, so that
 $$ q(u) > (1-\e)q(u^*),$$
 which combined with (\ref{neq1}) shows the first inequality of the lemma.  Repeating the same argument with $u$ and $u_*$ in place of  $u^*$ and $u$, respectively, we obtain $ q(u_*) > (1-\e)q(u)$, hence the second inequality of the lemma.  \qed 
\v2
We may analogously define
$$q_{K,U}^W(x) = P[\,W^x \; \mbox{exits from $U$ through }\; K\,].$$
%Cor2.3
\begin{Cor}\label{YY4.2}\, Let $\eta =\eta(\e,1)$ be as in Proposition \ref{YY4.1} and $u\in V(U)$ be such that 
$q_{K,U}^W(u)>0$.  Suppose that for any $\a>0$,  there exists positive constants $R_0>1$  and $1<\la \leq {\rm dist}(u,\partial U)$ such that  if $\rho\geq R_0$, then   $q_{K,U}^W(u)>0$ and   %for all $u\in V(\D)$ with ${\rm dist}(u,\partial U)>\e$,
\beqn\label{neq2}
1-\a\leq q_{K,U}^{(0)}(u; \eta\la)/ q_{K,U}^W(u)  \leq 1+\a.
\eeqn
Then  for any $\a>1$,  there exists $R'$ such that for $\rho\geq R'$, %and for all $u\in V(\D)$ with ${\rm dist}(u,\partial U)>\e$,
\beqn\label{neq3}
1-\a \leq q_{K,U}(u)/q_{K,U}^W(u) \leq 1+\a.
\eeqn
Here if $R_0$ is independent of  $u$,  $U$ and  $K$ (when these vary in any fashion), then so is  $R$.
 % Each of the first and second inequalities of (\ref{neq2}) implies the corresponding one of (\ref{neq3}). 
If  for some constant  $c$
\beqn\label{neq21}
q_{K,U}^{(0)}(u; \eta \la)/ q_{K,U}^W(u) < c  
\eeqn
in place of  (\ref{neq2}), then $q_{K,U}(u)/q_{K,U}^W(u) <2c$ in place of  (\ref{neq3}). 
\end{Cor}

{\sc Remark 1.} \,   The shifted walk $S^{0}\circ\th_{\sigma(\Ga_{\eta\la}(u))} $ behaves like Brownian motion $W^u$ (under  scaling by $1/\rho$) as long as they are kept  away from the boundary of $U$ with  some sufficient  distance.  In order to ensure condition  (\ref{neq2})  or something like that we need some condition for the pair $K$ and $U$. If $U$ is nice, such a  condition  will be satisfied for any  $K$.
   As is discussed in the next section we are concerned with a conformal map $\fa_D$  from a simply connected domain $D$ onto $\D$. For each $0<\de<1$ fixed,  any nice domain $A$ in  $(1-\de)\D$, $U= \fa_D^{-1}(A)$ will be also nice (cf. Lemma \ref{Fact})). 
   
  %   \2.4
  \v2\n
 {\sc 2.4 \; Supplement to the proof of Proposition \ref{YY4.1}.}
  \v2 
    The proof of Proposition 4.1 of \cite{YY}  is based on the following fact:  If  $\sigma_\Ga$ (or $\sigma(\Ga)$) denotes the first epoch  the walk $S^{h,0}$ enters into a set $\Ga$.  Let  $a\in \D$ and define events $\La(r)$, $r>0$   by
 $$\La(r)= [ \bar S^{h,0}\circ \th_{\sigma(\Gamma_{r/20}(a))} \; \mbox{{\rm encompasses}} \; \Gamma_{ r/20}(a) \; \mbox{{\rm in}}\; \Gamma_{r}(a)],$$
 where  $ \bar S^{h,0}\circ \th_\sigma$ denotes  the curve that linearly interpolates   
 $S^{h,0}\circ \th_\sigma$.
Then,  for any $r>0$, there exists $h_0 =h_0(r)>0$ 
such that for  $0<h<h_0$, 
\[
P\Big[\mbox{both $\La(r)$ and $\La(20r)$}\; occur   \, \Big|\, \sigma_ {\Gamma_{r/20}(a)} 
< \tau_{(1+r)\D} \Big] >c
\] 
with some universal constant $c>0$.
What we need to ensure is the independence of $h_0$ from $a\in \D$. To this end it is irrelevant whether  we consider the event  $\La(r)$ or  $\La(r) \cap \La(20r)$, and  for simplicity 
we  verify the following

%Lem2.4
\begin{Lem}\label{lem1} \; For any $\e>0$ there exists a positive number $\eta=\eta(\e)<1/2$ such that for any positive number $\la<1/2$  one can choose $h_0=h_0(\e,\la)>0$  so that if  $h<h_0$ and $z \in (1-\la)\D$, then 
$$
P\Big[ S^{h,0}\circ \th_{\sigma(\Gamma_{\eta \la}(a))} \; \mbox{{\rm encompasses}} \; \Gamma_{\eta \la}(a) \; \mbox{{\rm in}}\; 
\Gamma_\la(a)  \, \Big|\, \sigma_ {\Gamma_{\eta \la}(a)} < \tau_\D \Big]  >1-\e.
$$
%(Note that $\Gamma_\la(w)\subset D$ due to the premiss $\fa_D(w)\in (1-\la)\D$ so that $h \Gamma_\la(w)\subset \D$ due to $D\in \D_h$.)
\end{Lem}
\v2\n
\pf\;  For each  $s\in (1/2,1)$ let   $A_s$ and $B_s$ denote the events defined by
\[
A_s = \mbox{[Brownian motion  $W^0$   hits  $\Gamma_{ s\eta \la}(a)$ before exiting $s{\D}$]},
\]
\[
B_s = [W^0\circ \th_{\sigma^W(\Ga_{\eta\la}(a))} \; \mbox{ encompasses\; $\Gamma_{\eta \la/s}(a) $ in} \;\Gamma_{s \la}(a)],
\]
where $\sigma^W_\Ga$ denotes the first hitting time of $\Ga$ by $W^0$,
and put  $p(s) = P[A_s]$ and $q(s)=P[B_s|A_s]$. Then $q(1)\geq  1-\e/3$ if $\eta=\eta(\e)$ is chosen small enough, and $p(s)/p(1/s)\uparrow 1$ and  $q(s) \to q(1)$ as $s\uparrow 1$. Fix $s<1$ so that $p(s)/p(1/s)  >1-\e/4$ and $q(s)\geq 1-\e/4,$ 
Noting that the boundaries of  both events $A_s$ and $B_s$ are null  we apply 
 the assumed invariance principle (H) to see that    
  our random walk $S^0$ and the Brownian motion $W^0$ can  be  both  defined on the same probability space so that  if the event $C_h$ is defined by
$$C_h=[d_{\D}(\bar S^{h,0},  W^0) < (1-s)\eta \la],$$
then for some $h_0=h_0(\e,\la)$
\beqn\label{ABC}
P[C_h]\geq 1-\frac14 \e P[A_s\cap B_s] \qquad 0<h<h_0.
\eeqn
 In the definitions of  $A_1$ and $B_1$ replace  $W^0$  by  $ S^{h,0}$ and let    $A^{RW}_1$ and $B^{RW}_1$ be the corresponding events. Then, from (\ref{ABC}) it follows that
$$P[B^{RW}_1\cap A^{RW}_1] \geq P[B_s\cap A_s\cap C_h] \geq (1-\e/4)P[B_s\cap A_s]  =  (1-\e/4)q(s)p(s)$$
and 
$$P[A^{RW}_1] \leq P[A^{RW}_1\cap C_h] + (1-P[C_h]) \leq  P[A_{1/s}] + 4^{-1}\e P[A_s]\leq (1+\e/4)p(1/s); $$
hence
 the probability in question is bounded from below by  $(1+\e/4)^{-1}(1+\e/4)^3\geq 1-\e$ as required. \qed.

%%\section{ Preliminary facts of geometric nature}  
\section{ Preliminary lemmas of geometric nature }  

 Let $D$ be a simply connected  domain of the complex plane and  $\hat o$ and $v_0$ be  fixed points of $D$ and $V(\partial D)$, respectively.  By the Riemann mapping theorem there exists   a (unique) conformal  map of  $\fa_D$ onto $\D$ such that  
 \beqn\label{eq4.1}
 \fa_D(\hat o)=0\quad \mbox{ and} \quad  \fa_D(v_0) =1.
 \eeqn
% where the second equality is understood to mean that whenever  $x_n\in D$ converges to $v_0$, then  $\fa_D(x_n) \to 1$.
  To be  precise  $v_0$ must be understood to be  a prime end of  $D$: otherwise $v_0$ may correspond multiple points of  the unit circle $\partial U$, hence multiple  $\fa_D$'s, and  in such a case it is understood that any one of them is selected. (For instance if $D$ is an upper half plane  with a slit $[0,i]$,   every point  $si, 0\leq s<1$, should be counted twice.   It is known that any conformal map of $D$ onto $\D$ naturally induces one to one correspondence between the set of prime ends of $D$ and $\partial \D$ \cite{Pm}, in particular   if $D$ is a Jordan domain, the prime ends are identified with the boundary points. For more details see Section 5.4; until there we shall not encounter any serious problem that  necessitates to use the concept of  prime ends.)
  
In this section we collect certain simple geometric relations between   the subsets of  $D$ and  their images  $\subset \D$  by $\fa_D$. Although  the planar graph $G$ is irrelevant to the analysis of this section, 
the results obtained have consequences on  the random walk on it  which are  included  in this section as Corollaries \ref{cor3.1} and \ref{cor3.2}. 

For any non-empty set $A\subset \C$,  denote  by ${\rm dist}(A,B)$ the distance between $A$ and another  $ B \neq \emptyset$, by ${\rm diam}\,A$ the  diameter of $A$ and by $\mbox{in-rad}_{x} \,A$ the inner radius of $A$ with respect to $x\in A$: ${\rm dist}(A,B)= \inf\{|x-y|: x\in A, y\in B\}$, ${\rm diam}\, A= \sup\{|x-y|: x, y \in A\}$ and 
$\mbox{in-rad}_{x} \,A= {\rm dist}(x, \C\setminus A).$
We continue to denote 
by $\Ga_r(a)$  the open disc of radius $r$ centered at $a$ (as defined in  (\ref{disc0})).

% We  denote 
%by $\Ga_r(a)$  the open disc of radius $r$  with  center at $a\in \C$: 
%\beqn\label{disc0}
% \Gamma_r(a) =\{u\in \C: |u-a|<r \}.
% \eeqn

\v2\n
{\sc 3.1. \; Elementary bounds on distortion under $\fa_D$}   
\v2

Put
 $$ \rho_D =\mbox{in-rad}_{\hat o}\, D =  {\rm dist}(\hat o, \C\setminus D).$$
 A version of the  Koebe distortion  theorem \cite[Corollary 1.4]{Pm} says that if $x\in D$ and $\de= 1- |\fa_D(x)|$,
 \beqn\label{K1}
\de(2-\de)/4 \leq |\fa'_D(x)|{\rm dist}(x,\partial D)  \leq \de(2-\de).
 \eeqn
 Taking $x=0$, this gives  $ \rho_D  \leq 1/|\fa_D'(\hat o)| \leq 4\rho_D$ and, employing  another form of the  Koebe distortion  theorem \cite[Theorem 1.3]{Pm}, we obtain that   for $0<r\leq 1$, 
 $$ \fa_D^{-1}(r\D) -\hat o \supset \frac{r}{ (1+r)^{2}\fa'(\hat o)} \D\supset \frac1{4} r\rho_D \D;$$
  similarly for $0<r< 1$,  $\fa_D^{-1}(r\D) - \hat o  \subset [4r\rho_D/(1-r)^2]\D$. 
  
 %We shall give several estimates  independent of  $D\in \D$ of  some probabilities, for which  the following fact  is fundamental.
 \v2\n
 %%Lemma \ref{Fact}  
 \begin{Lem} \label{Fact}
  There exists a positive increasing function $\k(\ell)$ on the interval   $ 0<\ell <1$ such that 
 if $K\subset   D$ is a compact connected set  and ${\rm diam}\, \fa_D(K) \geq\ell$, then ${\rm diam}\,  K \geq  \k(\ell) \rho_D,$ entailing that  if a line  segment 
 $[x,y]$ is contained in $\overline D$ and $|x-y| < \k(\ell)\rho_D$, then $|\fa_D(x)-\fa_D(y)|\leq \ell$.   It in particular follows that if  $z\in D$ and $1- |\fa_D(z)| \geq \de$, then 
 $${\rm dist}(z,\partial D)\geq  \k(\de)\rho_D  \quad \mbox{and} \quad \fa'_D(z)| \leq [(2-\de)\de/\k(\de)]/ \rho_D.$$ 
 \end{Lem}
 \v2\n
 \pf\; 
The distortion theorem says that
 $| (\fa_D^{-1})'(z)\fa_D'(\hat o) |  \geq (1-|z|)/ (1+|z|)^3\geq \de /8$ for $|z|<1-\de $, which combined with the inequality  $\rho_D \leq 1/|\fa_D'(\hat o)|$  shows the first bound  of the lemma  (with $\k(\ell)=\ell \de/8$)  if $\fa_D(K) \subset (1-\de)\D$.
Consider the case $\fa_D(K) \not \subset (1-\de)\D$. We  suppose $\fa_D(K) \not \subset \frac12\D$ for definiteness.
Let   $r_*$ be the infimum of  $r$ such that $\fa_D(K) \subset r\D$ and put $D_* = \fa_D^{-1}(r_*\D)$.
 Then $r^*\geq 1/2$, hence
  $\mbox{in-rad}_{\hat o}(D_{*})\geq \frac14 r_*\rho_D \geq  \frac18 \rho_D$ and   from  the Beurling estimate it follows that
the harmonic measure of $K$   in   $D_*\setminus K$  from $\hat o$  is bounded above by $c'\sqrt{\rho_D^{-1}{\rm diam}\,  K }$ with a universal constant $c'>0$ (\cite[Corollary  3.78 (the second formula)]{L}), whereas
 the same harmonic measure is bounded below by a positive multiple of $ {\rm diam}\,  \fa_D(K) $  owing to the conformal invariance of the harmonic measure, since  $\fa_D(K)/r_* $ is a connected subset of the closed unit disc $\overline \D$.
 This entails that 
 $$ {\rm diam}\,  \fa_D(K) \leq c''\sqrt{\rho_D^{-1}{\rm diam}\,  K },$$
 which gives the required inequality with $\k(\ell)= (\ell /c'')^2$.

\v2\n
{\sc Remark 2.} \; We may take $\k(\de)= c\de ^2$ with a universal constant $c>0$ as is indicated  in the proof, but we do not need it in this paper.

\v2\n  
{\sc 3.2.\; Domains $ B_r(x), U_\de$ and $ Q(x)$.} 

\v2
For $a\in D$ and $r>0$ let $ B_r(a)= \fa_D^{-1}(\Ga_r(a))$, namely
$$ B_r(a)=\{u\in D: |\fa_D(u)-\fa_D(a)|<r \}.$$ 
For $\de>0$ put
 $$U_\de  =U_\de ^{D, v_0} = \{x\in D: {\rm dist}(\fa_D(x),\partial \D)<\de  \}$$
 and  
 $$U_{r,\de }= U_{\de } \cap B_{r+\de}(v_0).$$ 
  %In below we write   $d^*_D(x, \partial D)= {\rm dist}(\fa_D(x),\partial \D)$.  
 By  Lemma \ref{Fact}  
 \beqn\label{fact0}
 {\rm dist}(\partial U_\de, \partial D) \geq \k(\de)\rho_D.
 \eeqn
    
  We adapt a method found  in \cite{LSW}. 
  Given  $\de>0$ and  $x \in D\cap \partial U_{\de}$ and let $z^*=z^*(x,\de)$ be a point of $\partial D$ closest to $x$ and set
$r= {\rm dist} (x,\partial D)= |z^*-x|$. 
% let $C_{\k}=C_{\k}(x,\de)$ be the connected component of $D\cap \partial \Ga_{\k r}(z)$,
 For $0 < \k<2$, let  
\beq
 Q_{\k}= Q_{\k}(x,D) &=& \mbox{  the connected component of}\,\, \Ga_{\k r}(z^*)\cap D  \\
&& \mbox{ which contains a point of the segment $[z^*, x]$.}
 \eeq  
We write $Q(x)$ for $Q_1=Q_1(x,D)$. 
 The following result (as well as   its proof)   is a  simple modification of that found in the proof of Lemma 5.4 %(asserting (\ref{hm})) 
 in \cite{LSW}. % It is  fundamental in the arguments in this section.
  (The modification, although not substantial  at all,   make simpler and clearer  the arguments developed  later.)
  
 %%%Lem{LSW} 
 \begin{Lem}\label{LSW} \;   Let  $Q(x)= Q_{1}(x,D)$  be defined as above,  $\om$ the component of $\partial Q(x)\cap D$ containing  $x$ and  $D(x)$ the component of $D\setminus \om$ that does not contain $\hat o$.
There exists a universal   constant  $m > 1$  such that if    $0<\de <1/m$ and $x \in D\cap \partial U_{\de}$, then   $\om\cup \Ga_{r/4}(x)$ and $\partial B_{m\de}(x)\cap D$   are disjoint;  in particular    only the following  two alternatives   are possible:  
\beq
 (1)\quad Q(x)\subset  D(x);\qquad (2) \quad  Q(x) \subset D\setminus D(x);
  \eeq
 %it holds that     \beqn\label{hm0}
%B_{m\de}(x)  \supset  \left\{ \begin{array} {ll} Q(x)\cup \Ga_{ r/4}(x)   &\mbox{ in  the case {\rm (1)},} \\
%   D(x)\cup \Ga_{ r/4}(x)  \quad &  \mbox{ in the case {\rm  (2);}}   \end{array}   \right.\eeqn
 in either case
   \beqn\label{hm0}
D(x)\cup \Ga_{ r/4}(x) \subset B_{m\de}(x). 
  \eeqn
  %and for some universal constant $c>0$,
 %\beqn\label{hm}
 %H^{BM}(x, \partial D; B_{m\de}(x)) > c,
 %\eeqn
 %where    $H^{BM}(y,A; \Om)$ denotes the harmonic measure of  $A\subset \C\setminus \Om$ in a domain  $\Om$ from  $y\in \Om$.
%(The inclusion    $Q(x)\cup \Ga_{ r/4}(x) \subset B_{m\de}(x)$ holds always  in the case (1); it may and may not be true in the case (2).)
\end{Lem}
\v2\n
\pf \; %The proof  is similar to the proof of Lemma 5.4 of \cite{LSW}. 
 By the Koebe $1/4$ theorem  we have i-rad$_x B_\de(x) \geq \frac14 \de(2-\de)/|\fa_D'(x)|$, which combined with (\ref{K1}) gives the inequality  $r\leq \frac14 \mbox{in-rad}_x B_\de(x)$ so that for all $m>1$,
 \beqn\label{2.5}
\overline{ \Ga_{ r/4}(x)} \subset  B_{m\de}(x).
 \eeqn
Let $A_+$ denote the annulus $\Ga_{\frac54 r}(z^*) \setminus \Ga_{ r}(z^*)=\{z: r\leq  |z^*- z|<\frac54 r\}$.  Observe the following two simple facts, the former is obvious from the definitions of $A_+$ and $Q(x)$;  the latter follows from the conformal invariance of harmonic measure.  

%
%\begin{figure}[t] 
%\begin{center}
%  \vspace*{-3.0cm} 
%   \hspace*{0em} 
%\includegraphics[width=11cm,height=7.5cm,clip]{top1.eps}
% \vspace*{-0.2cm} 
%\end{center}
%\caption{Figure I}
%\end{figure} 

%\begin{figure}[t] 
%\begin{center}%%%%%\centering
%\vspace{1mm}
% \hspace*{-1.0em} 
%\includegraphics[width=10cm,height=4.9cm,clip]{ushar7.eps}
%\vspace*{-2mm}
%\end{center}
%\end{figure} 

 (a) With a  probability greater than a positive universal constant the Brownian path  $W_t$  started at $x$ moves around  $Q(x)$ counterclockwise  in such a way that   for some  epoch $t_0>0$ the path $W[0,t_0]$ is contained  in $ A_+\cup \Ga_{r/4}(x)$ and as $t$ ranges over $[0,t_0]$ its argument   about  $z^*$, $\arg (W_t-z^*)$,   extends  through at least an interval  of length $11\pi/6$ while confined in $(\th-\pi/6 ,  \th + 11\pi/6)$ where $\th =\arg(x-z^*)$; by symmetry the same thing but in the opposite direction of rotation holds. %(See Figures 1 and 2 where  $x=\fa^{-1}(1-\de)$.)
\v2
  (b)  The probability that the Brownian motion started at $x$ leaves $B_{m\de}(x)$ before hitting $\partial D$ may be made arbitrarily small for  $\de <1/m$  by choosing  $m$ large.  
\v2
Taking (\ref{2.5}) into account  it follows from (a) and (b) above  that 
if $m$ is large so that  both  the probabilities in (a) are larger than that in  (b),  then $\partial B_{m\de}(x) \cap D$ cannot disconnect   $x$ from $\partial D$ in either of the two components  of  $Q_{5/4 } \setminus (Q(x) \cup \Ga_{ r/4}(x))$   since otherwise  any encompassing path in (a) must hits $\partial B_{m\de}(x)$ earlier than  $\partial D$ (entailing that  the event in (b) is contained in one of two events  in  (a)).

In the same way, but considering the annulus  $A_- = \Ga_{ r}(x) \setminus \Ga_{\frac34  r}(x)$ in place of $A_+$, we  see that for  $m$  large  enough   $\partial B_{m\de}(x)\setminus \partial D$ cannot disconnect  $x$ from $\partial D$ in either of the two components  of  $Q(x) \setminus (Q_{3/4} \cup \Ga_{ r/4}(x))$.   

Simple topological arguments then verify  that  $\partial B_{m\de}(x)\cap (\om\cup \Ga_{r/4}(x))=\emptyset$  and either (1) or (2) holds.   \qed

\v2
%Cor3.3
\begin{Cor}\label{cor3.1}  For any $\de_0$ and $M$ there exists  a constant  $R$  (depending only on $\de_0$ and $M$) such that  if $D\subset M\rho_D \D$ and $\rho_D\geq R$, then  for some  universal constant $c>0$,
$$
  P[ S^{ u} \;\mbox{exits  $B_{m\de}(u)$  through  $\partial D$}\,]  >c \qquad \mbox{if }\;\; u\in V(D\setminus U_{\de_0}),
$$
where  $\de= 1-|\fa_D(u)|$  and $m$ is the universal constant  appearing  in Lemma \ref{LSW}.
\end{Cor}
\v2\n\pf\,  We apply  Lemma \ref{LSW}  with   $u$ in place of $x$. Let  $A_-$ be the annulus defined in its proof. Then in the case (1) of Lemma \ref{LSW} the connected component of $(A_-\cap D)\cup \Ga_{r/4}(u)$ is contained in $B_{m\de}(u)$,  so that 
$$P[ S^{ u} \;\mbox{exits  $B_{m\de}(u)$  through  $\partial D$}\,]  \geq P[S^u\; \mbox{encompasses} \; \Ga_{3r/4}(z^*)  \; \mbox{in} \; A_-\cup\Ga_{r/4}].$$
Owing to   Corollary \ref{YY4.2} (and the remark after it) the right side is bounded below by  half the corresponding probability for Brownian motion, which is a universal constant. 
In the case (2), a similar reasoning shows the result. \qed

\v2
\n
{\sc 3.3. \;  Construction of a domain $U$. }

\v2
  We give a consequence of Lemma \ref{LSW} in a form that is actually used in our substantial application of Lemma \ref{LSW}.
Let  $m$ be the universal   constant  described in Lemma \ref{LSW}. 

\v2\n
%%%Lem\label{LSW2}
\begin{Lem}\label{LSW2} \;    For each  $0<\de<1/m$,    there exists a simply connected subdomain $U$    of $U_{m\de,\de}$ such that both of the sets
\beqn\label{C_C}
C:= \partial U\cap \partial D,\quad C^*:= \partial U\cap (D\cap \partial U_\de) \quad
\eeqn
are not empty and that  if
$$\Om = \C\setminus  (\partial U \setminus C),\quad$$
then for some $x^*\in C$ and some universal constant $c_*>0$,
\beqn\label{est_H}
h^{bm}(y, C^*_{in}; \Om) \geq c_* \quad \mbox{for}\quad y\in \Ga_{\k(\de)\rho_D/4}(x^*),
\eeqn
where $\k$ is the function given in Lemma \ref{Fact}, $C^*_{in}$ means the inner side of   $C^*$ relative to $U$, and  $h^{bm}(y,A;\Om)$ denotes the harmonic measure of  $A$ in  $\Om$ from $y$ (so that 
$$h^{bm}(y, C^*_{in};\Om)= P[ W^{y}(\tau_\Om) \in C^*\;\;\mbox{and}\;\;  \exists\, \e>0,  \forall\, s\in (0,\e), W^{y}(\tau_\Om-s)\in U]\;).$$
\end{Lem} 
\v2\n
\pf  \;
  Take  $\fa^{-1}_D(1-\de)$ for $x$  in  Lemma \ref{LSW},  which we are to apply.  Let $Q=Q(D,x)$,    
  $z^*=z^*(x,\de)$  and $r$ ($= |x-z^*|$) be as in Section 3.2 and $D(x)$ and $\om$ be as defined in Lemma \ref{LSW}. 
  Our construction of $U$ apply for either of  two alternatives  in Lemma \ref{LSW}, while $x^*$  is chosen  in a different fashion.

%In  Lemma \ref{LSW} we have two alternatives. First consider the case (1), i.e.,   $Q(x)\cup \Ga_{ r/4}(x) \subset B_{m\de}(x)$.  
     Put $U^0 = U_\de \setminus [ D(x)\cup \Ga_{ r/4}(x)] $ (recall  $ D(x)$ is the component of $D\setminus \om$ that does not contain  $\hat o$)  and   define $U$ by
\beqn\label{d_U}
U =  \mbox{the component of\; $ U_\de \setminus U^0$ \; that contacts  $\partial D$}.
\eeqn
%(See Figure 3.) 
$U_\de\cap\partial U$ consists of two disjoint arcs,  $C_\pm$ say,  of $\partial [D(x)\cup \Ga_{ r/4}(x)]$;
and  $C$  (given by  (\ref{C_C})) is  the  part of  $\partial D$  (in contact with $U$) cut  by the endpoints of $C_\pm$  that agree those of $\om$. (To be precise either of the endpoints of $\om$  do not always cut  $\partial D$ if the cluster set to the prime end associated with it is not  point, but this causes no problem in below.)  Similarly $C^*$ is an  arc cut from the inner boundary of $U_\de$ by the other endpoints of $C_\pm$. Clearly  $\fa_D(U) \subset \Ga_{m\de}(1)$, and we see that   $U\subset U_{m\de, \de}$. 

Now consider the case (1) and we take $x^*=z^*$, which is certainly in $C$.  
 If the intersection of $\Ga_{r/2}(z^*)\setminus Q(x)$ with  $U$ is empty, then ${\rm dist}(\partial \Om, z^*)\geq r/2$,  and  it is  plain to see (\ref{est_H}). 
If this intersection is  not empty, we modify $U$ by subtracting  the set $\overline{\Ga_{r}(z^*)\setminus Q(x)}$ from it and taking the component of the resulting set that intersects $\Ga_{r/4}(x)$ and define  $C^*$ by  (\ref{C_C}). Then we also have (\ref{est_H}).

In the case (2) 
we take as  $x^*$  the  point where the ray  that  issues   from $z^*$ and passes through $x$ first falls on  $\partial D$. %; if the ray does not intersect $\partial D$ we let  $x^*$ be either of the  endpoints of the component  of  $D\cap\partial \Ga_{R}(x)$ that intersects the ray, where  $R= (4r)\vee (\k(\de)\rho_D)$.  
Then, on observing that any  continuous
curve  from $x^*$  that reaches  $\partial \Om_{out}$  (the $z^*$ side of $\partial \Om$)   within  $\Om$ must travel a distance more than $\pi r/3$, we see  (\ref{est_H})  holds true.  The proof of Lemma \ref{LSW2} is finished. \qed

%\begin{figure}[t]
%\begin{center}
%\vspace*{-12mm}
% \hspace*{0.0em} 
%\includegraphics[width=10cm,height=6.4cm,clip]{ushar3.eps}
%\vspace*{-2mm}
%\end{center}
%\end{figure} 

\v2
Let $H^{RW}(u, C^*_{{\rm in}}; \Om)$ denote the harmonic measure of $C^*_{{\rm in}}$ for the walk $S^u$
in $\Om$.  (For the present purpose the detailed definition is irrelevant  and any reasonable one may be adopted.) Then we have
\v2
%Cor3.5
\begin{Cor}\label{cor3.2} For any $\de_0>0$ and $M\geq 1$ there exists  a constant  $R$  (depending only on $\de_0$ and $M$) such that   if $D\subset M\rho_D \D$ and $\rho_D\geq R$, then  for some  universal constant $c>0$,
$$
 H^{RW}(u, C^*_{{\rm in}}; \Om) >c\quad \;\mbox{for}\quad\; u\in \Ga_{\k(\de)\rho_D/8}(x^*)\cap V,
$$
where $m$ is the universal constant  appearing  in Lemma \ref{LSW}.
\end{Cor}
\v2\n\pf\, Comparing the harmonic measure  on the left side of  (\ref{est_H})  with the corresponding harmonic measure
for the random walk $S^{0}\circ \th_{\sigma(\Ga_{\k(\de)\rho_D/8}(y))}$ we infer that the latter, which may be written as 
$q^{(0)}_{C^*_{{\rm in}},\Om}(u; \k(\de)\rho_D/8)$, is larger than $2c$  if $\rho_D$ is large enough,   and an application of Lemma \ref{Nlem1} concludes the proof. \qed

 %P[ \bar S^{u}(\tau_\Om) \in C^*\;\;\mbox{and}\;\;  \exists\, \e>0,  \forall\, s\in (0,\e), W^{y}(\tau_\Om-s)\in U]\;). 
 %%%%%%%%%%%%%%%%
%%%%%%%%%%%%%
%%%%%%%%%%%%

%\section{Estimates of hitting distributions of random walks }
\section{Estimates of hitting distributions of random walks}

Let $G= (V,E)$ be the planar graph and $(S^v_n)_{n=0}^\infty$ the random walk on  $V$ starting at $v$ described in Section 2.
A  bounded domain $D$ is called a grid domain if its boundary consists of edges of the graph $G$. Define
 $$\cD =\{D: D \; \mbox{is  a simply connected and bounded grid domain}  \}.$$
In this section we shall give several estimates of harmonic measures of the walks $S^v$,   $v\in D$ for various subdomains of $D\in \cD$ that are defined by means of a conformal $\fa_D: D\mapsto \D$.  What causes the problem is that the map $\fa_D$ may distort the metric property of $D$ unrestrictedly when the points approach the boundary of $D$, so that the direct application of the  invariance principle---such   as that stated in Corollary \ref{YY4.2}---would be  impossible since we must verify   condition (\ref{neq2}) or (\ref{neq21}).

 In what follows,  as in the preceding section,   a point $\hat o$ of  $D\in\cD$ is suitably chosen as a \lq reference point'  which together with $v_0\in \partial D$ uniquely determines   the conformal  map  $\fa_D: D\mapsto \D$ via   the condition (\ref{eq4.1}). With this being taken into account define
 $$\cD_{R,M}= \{ (\hat o,D): \hat o\in D\in \cD, v_0\in V(\partial D),  \rho_D >R, D\subset M\rho_D\,\D\},$$
 where   $\rho_D= \mbox{in-rad}_{\hat o}\, D$ as in Section 3.
 For sake of  brevity we write $D\in \cD_{R,M}$ instead of $(\hat o, D)\in \cD_{R,M}$ with the understanding that  $\hat o$ is assigned  to $D$  in  some  way (cf. the beginning of Section 5.2 and   (\ref{p(t)}) for the choice of $\hat o$). We shall be interested in a lower bound of $\rho_D$ and apart from it  no significance will be attached  to a particular choice of  $\hat o$ for the discussion made in the present section. A boundary point  $v_0\in V(\partial D)$ is also supposed to be assigned to $D\in \cD$  to determine $\fa_D$ uniquely, but our analysis will be carried out so as to be independent of it.  The dependence on $M$ is needed because   the   bounds   in Proposition   \ref{YY4.1} or in (\ref{neq2}) (with a nice $U$) are not ensured by (H) if $u$ is indefinitely far from the origin.  If    these estimates are  valid without the restriction $|u|<M$, then it is unnecessary  to impose the boundedness  of $D$ indicated by  $M$. In any case we shall do not  take much care of the restriction imposed by $M$.
 
  Throughout this and the next section  we continue to use the notations $U_\de, U_{r,\de},$ $ B_r(x)$ and $\Ga_r(x)$ (the first three  are defined in Section 3.2 and the last by (\ref{disc0})); also suppose the condition (H) introduced  in Section 2 to be valid.  The subscript $D$ is dropped from $\fa_D$ in the proofs, if  $D$ is clear from the context. In addition we bring in 
\beqn\label{bar_tau}
\bar \tau_U = \sup\{ t > \tau_U -1 : \bar S^v_t \notin U\}.
\eeqn
For instance the statement that  $S^v$ exits $U_{r,\de}$ through  $\partial U_\de\cap U_\de$ is expressed as  $\bar S^v(\bar \tau_{U_{r,\de}}) \in U_\de$, while $S^v(\tau_{U_{r,\de}})$ is possibly  in  $U_{r,\de}$. Note that  
the expressions $S^v(\tau_{U_{r,\de}}) \notin D$, $S^v(\tau_{U_{r,\de}}) \in \partial D$ and $\bar S^v(\bar \tau_{U_{r,\de}}) \in \partial D$  all mean the thing if $D$ is a grid domain.
\v2\v2\n
{\sc  4.1. Simple properties  of the planar graph}
\v2

 Here we state some elementary  results that  follow from Hypothesis (H);
\v2

(1)\; $\max\{|u-v|: (u,v)\in E, u, v\in V( r\D)\} =o(r)$\; as \;$r\to\infty$.
\v2
(2)\; $ \sup\{{\rm dist}(z, V): z\in r\D\} =o(r)$\; as \;$r\to\infty$.
\v2
(3)\; ${\displaystyle \sup_{D\in \cD_{R,M} }\sup_{\,(u, v)\in E,  u\in \overline D, v\in D} |\fa_D(u)-\fa_D(v)| \to 0  }$\; as \;$R\to \infty$ (for each $M$). 
\v2
(4)\; For any $\de>0, \e>0$ and $M>1$, there exists $R>1$ such that if $D\in \cD_{R,M}$, $|\fa_D(x)| \leq 1-\de$, then 
 there is  a path of $G$   of diameter  less than $\e\k(\de)\rho_D$  that encircles $x$ in  $D$ and   $B_\e(x)$ contains a vertex of $V(D)$. ($\k$ is the function specified in Lemma \ref{Fact}.)
\v2

(1) and (2)  are readily verified directly from (H). For  (3)  use Lemma \ref{Fact}  on  noting that if $(u, v)\in E,  u\in \overline D$ and  $v\in D$, then $[u,v]\subset \overline  D$ in view of planarity of $G$ and the definition of grid domain. (4) also follows from Lemma \ref{Fact} together with (H). %combined with Corollary \ref{YY4.2}.

The  facts listed above, which are  easy to grasp, will be applied without  explicitly mentioning  of  their use in most cases in the later discussions.  From them    it follows  that given  $M>1$, $\eta\in (0,1/10)$ and $\de>0$, we can choose $R>1$ large enough  that for every  $y\in \partial D$ the  subgraph $G(B_{\de}(y) \setminus U_{\eta \de})$ is so  spatial  as to contain a fine network of paths: e.g.,  the annulus $B_{\de}(y)\setminus B_{\de/2}(y)$  contains a path of  $G$ that connects the two disconnected parts of 
$\overline {B_{\de}(y)\setminus B_{\de/2}(y)} \cap U_{\eta \de}$. It will be tacitly supposed  that $R$ is  large enough according to the arguments developed in below.  

  %$R$ is tacitly supposed to be chosen large enough.

\v2\v2\n
{\sc  4.2. Starting near the boundary (unconditional case) }
\v2

The walk starting near the boundary is relevant to our analysis. In this subsection we verify that  the probability of such a walk escaping immediate absorption to the boundary is small and use this fact to derive a result on the harmonic measure of the walk in $D\in \cD_{R,M}$.  In the next subsection we consider the behavior of the walk conditioned to escape immediate absorption.

  The next lemma is a slight improvement of Proposition 4.5  of \cite{YY}  (for the present setting) and the corresponding one verified in the proof of Lemma 5.4 of \cite{LSW} (for the simple random walk on the square lattice).  The proof is similar except for  our use of Lemma \ref{LSW}.  (An extended form is found  in the proof of Corollary \ref{prop5.2}.)
  
  %Lem4.1
  \begin{Lem}\label{lem2.20}  For any $\de_0>0$  and  $M>1$  there exists $R$ such that  for  $\de\geq \de_0$, $D\in \cD_{R,M}$  and $u \in V(U_{ \de})$,
  \beqn\label{eqlem5.10}
   P[ S^{ u}(\tau_{ B_{m\de}(u)}) \in  \partial D] > c,
   \eeqn
   where $m$ and $c$ are  the same universal constants specified in Lemma \ref{LSW} and Corollary \ref{cor3.1}, respectively.
  \end{Lem}
  \v2\n
  \pf\;  Owing to   Corollary \ref{cor3.1}    there exists  a constant  $R$  (depending only on $\de_0$ and $M$) such that  for $\de\geq \de_0/3$,
\beqn\label{eq5_11}
  P[ S^{ u}(\tau_{ B_{m\de}(u)}) \in  \partial D]  >c, \quad \mbox{if }\quad u\in U_{\de}\setminus U_{\de/5}.
  \eeqn
%showing  (\ref{eqlem5.1}) with  $\e= 1-c/2<1$. %Here $H^{RW}_{K}(u, A)$ denotes the harmonic measure of the random walk on $G$. 
 
 Let $u\in U_{\de /5}$ and  put  $p(v)=P[S_{\tau(B_{m\de}(u))}^v \in \partial D]$ for $v\in B_{m\de}(u)$. Since $p(v)$ is harmonic in $B_{m\de}(u)$, a maximum principle shows that there exists a path $\ga$   of
  $G( B_{m\de}(u))$ that connects $u$ with  $\partial B_{m\de}(v)\setminus \partial D$ such that $p(v)\leq p(u)$ for $v\in \ga$.  (Here the phrase that  $\ga=(\ga_k)_{k=0}^n$   \lq connects'  $u\in V$ with  a set  $A\subset \C$ means  that $\ga_0=u$ and   $(\ga_{n-1}, \ga_n] \cap A \neq \emptyset$.)  If   $\ga$ intersects the curve $\partial U_{\de/4}\cap D$ (in $B_{m\de}(u)$), 
then $p(u) > c$ according to what we have shown in the preceding paragraph. 

 Consider the case when  $\ga$ does not, and   suppose $\fa_D(u)$ is real and positive for  convenience of description.  
Put  $U^\pm_\de= \fa^{-1}(\D^\pm) \cap U_\de$, where $\D^\pm = \{z\in \D: \pm \Im z >0\}$.
By replacing   $R$ by a larger one if necessary we can then find  two vertices $u_{+}$ and  $u_{-}$  such that  $$u_\pm \in U^\pm_{\de/3}\setminus  \overline{B_{m\de/3}(u)\cup U_{\de/4}} \quad \mbox{and}\quad B_{m\de/3}(u_{\pm}) \subset   B_{m\de}(u).$$
 Now let $\de\geq \de_0$ so that we may apply (\ref{eq5_11}) with $u_\pm$ and $\de/3$ in place of $u$  and  $\de$, respectively.  Then, we infer  in the same way as  above that  there are paths
   $\ga_\pm$ such  that $\ga_\pm$ connects  a vertex  $u_\pm$  with $\partial D$  in $ B_{\de/3}(u_\pm)$ (respectively) and $p(v) > c$ for $v\in \ga_+\cup \ga_-$.
If $\ga$ does not intersect the curve $\partial U_{\eta_0\de/4}\cap D$, then $\ga$ must cross either $\ga_+$ or $\ga_-$, hence 
   $p(u) >c$. The proof of Lemma \ref{lem2.20} is  complete. \qed
   
     %Cor4.2
 \begin{Cor}\label{cor2.2}  There exists a universal constant $c_1<1$ such that for any $\de_0>0$  and  $M>1$ there exists a constant $R$ such that for $\de_0\leq \de<1/4$, $D\in \cD_{R,M}$  and $u\in  V(U_{ \de})$ with $|  \Im\fa_D(u)|< \de/4$, 
  $$ P[ \bar S^{ u}(\bar\tau_{U_{\de, \de}}) \in  U_\de] < c_1.$$
  (See (\ref{bar_tau}) for $\bar \tau_U$.) 
 \end{Cor}
   \v2\n
  \pf\;   Given $\de_0>0$, take  $\de\geq \de_0$ and  $u\in  U_{ \de}$.  By Lemma \ref{lem2.20}   $ P[\bar S^{ u}(\bar\tau_{U_{\de, \de}}) \notin  U_\de]  > c$ whenever   $u\in  U_{\de/2m}$ and $|  \Im\fa(u)|< \de/2$.  The general case of  $|  \Im\fa(u)|< \de/4$ is reduced to what is just verified.  Indeed, on  putting $U= U_{\de/2, \de}\setminus U_{\de/2m}$, with the help of Lemma \ref{Fact}   an application of   Corollary \ref{YY4.2} shows that  if  $|  \Im\fa(u)|< \de/4$, then $P[\bar S^{ u}(\bar \tau_{U}) \notin  U_\de\setminus U_{\de/2m}]> c'$ with some universal constant $c'>0$, provided  that $R$ is large enough, so that the required inequality ensues at least with $c_1 = 1- c \,c'$.
\qed
  \v2
  
  %Cor4.3
 \begin{Cor}\label{lem2.2}
 For  any $\e>0$, $\a>0$  and  $M>1$  there exists $R$ and $\de>0$ such that  for  all   $D\in \cD_{R,M}$  and $u \in V(U_{\de})$,
  \beqn\label{eqlem5.1}
   P[ S^{ u}(\tau_{ B_{\a}(u)}) \in  D] < \e.
   \eeqn
 \end{Cor}
  \v2\n
  \pf\;  Taking $N$ so large  that   $\e <(1-c)^N$, we have only to  apply Lemma \ref{lem2.20} repeatedly  at most  $N$ times  by starting  with   $\de =\de_0= m^{-N}\a$ to  arrive at the inequality of the lemma. \qed
\v2  
  By virtue of the bound (\ref{eqlem5.1})  we  can control the probability  of  the random walk  badly behaving 
  near the boundary, and  an application of Corollary \ref{YY4.2}  leads to    the following corollary. 
       %Cor4.4
 \begin{Cor}\label{cor5.3}  For any $\e>0$, $r>0$ and $M>1$ there exists $R>1$ such that
 if $I\subset \partial \D$ is an arc  of length $2r$ centered at 1 and $J= \fa_D^{-1}(I)$, then  for $u\in V(D\setminus U_\e)$,
  \beqn\label{eqcor5.3}
   (1-\e)  P[ W^{ u}(\tau^W_{D}) \in  J ]  \leq P[ S^{ u}(\tau_{D}) \in  J ] \leq (1+\e)P[ W^{ u}(\tau^W_{D}) \in  J ] ,
   \eeqn
 where $W^u$ denotes a planar Brownian motion started at $u$ and $\tau^W_B$ its first exit time from $B$.
 \end{Cor}  
\v2\n
\pf\;  Let $0<\eta <1$ and for  $u\in V(D \setminus U_\e)$  and $0<\de < (r \wedge\e)/2$, denote by $A=A_{\de, \eta,u}$  the event that 
the walk enters into  $U_\de$ substantially    through $\partial U_\de\cap \partial  U_{(1+\eta)r,\de}$:
\beqn\label{A}
A=\{ \, S^u(\tau_{D\setminus U_\de})\in U_{(1+\eta)r, 2\de}\, \} \nonumber
\eeqn
and make decomposition
\beqn\label{dcm_A}
 P[ \, S^u(\tau_{D})\in J]  =  P( \{ S^u(\tau_{D})\in J\} \cap A)   +   P( \{ S^u(\tau_{D})\in J\} \setminus A). 
\eeqn
 According to Corollary \ref{eqlem5.1}  for any $\e_1>0$ there exists $\de_1 = \de_1(\eta r, \e_1)>0$  and $R=R(\eta r,\e_1)$ such that  $P[ \, S^w(\tau_{B_{\eta r}(w)})\in D]  <   \e_1 $ for $w\in U_{\de}$,  $\de\leq \de_1$, implying  that the second probability  of the decomposition  is less than   $\e_1$. On taking $\e_1=\frac12 \e\inf_{w\in V(D\setminus U_\e)} P[ W^{ w}(\tau^W_{D}) \in  J ] $ this yields that 
\beqn\label{2ndtrm}
 P( \{ S^u(\tau_{D})\in J\} \setminus A) \leq   \frac12 \e P[ W^{ u}(\tau^W_{D}) \in  J ]. 
 \eeqn
On the other hand,  applying  Corollary \ref{YY4.2}, we deduce that for each  $\eta$ and  $\de \leq \de_1$, we can choose    $R$  large  so that
$$P(A)  \leq  (1+{\textstyle \frac14}\e) P[ \, W^{u}(\tau^W_{D\setminus U_\de})\in  \partial U_{(1+\eta)r,\de}].$$
By the conformal invariance of Brownian hitting probability  the constants $\de$ and $\eta$  may have been chosen small (independently of $R$) so that the right side above is at most $(1+{\textstyle \frac12}\e) P[ W  ^{u}(\tau^W_D)\in J]$.
As a consequence, we may assert that for $R$  large enough,
%the  first probability on the right-hand side of (\ref{dcm_A})   is dominated  by
$$ P(A) \leq  (1+{\textstyle \frac12}\e) P[ W  ^{u}(\tau^W_D)\in J].$$
 This together with (\ref{2ndtrm}) 
concludes the upper bound of (\ref{eqcor5.3}).
 % since $\e$ is assumed so small that $(1-\frac14 \e)^{-1}(1+\frac12 \e)<1+\e$.

In a similar way, putting  $q_\eta(u) =P[ S^u(\tau_{D\setminus U_\de})\in U_{(1-\eta)r, 2\de}\,] $,  we infer  on the one hand
$q_\eta(u) \leq P[ S^{ u}(\tau_{D}) \in  J ]  +  {\textstyle \frac12}\e q_\eta(u)$, and on the other hand 
 $q_\eta(u)  \geq P[ W  ^{u}(\tau^W_D)\in J] (1-  {\textstyle \frac12}\e)$, yielding the lower bound in (\ref{eqcor5.3}).
 The proof is complete. \qed
 \v2
 
 A simple modification of  the  proof above shows extensions  of  (\ref{eqcor5.3}) that provide  for a certain class of events $B$ the upper and lower bounds of  $P(B\cap \{S^{ u}(\tau_{D}) \in  J\})$ in terms of the corresponding probability for Brownian motion. Among them we shall need the following 
one.  

%cor4.5
\begin{Cor} \label{ncor4.5} \,  Let $J$ be as in Corollary \ref{cor5.3}. Let   $u\in D\setminus  U_{2\de}$.  Let $\Ga$ be a set of paths of $G(D)$ and $\Ga'$ a measurable set of  continuous curves  in $D$ such that for  any  $\de>0$ small enough,  one can choose $R$ so that  for some constant $c>0$,
\beqn\label{Ga}
 P[S^u[0,\tau_{D\setminus U_{ \de}}] \in \Ga \; \;\mbox{and}\;\; S^u(\tau_{D\setminus U_\de})\in U_{\frac12 r,2\de}\,] 
 > c P[ W^u[0,\tau_{D\setminus U_{ \de}}]\in \Ga'  \; \;\mbox{and}\;\; W^{u}(\tau^W_D)\in J ].
\eeqn
  Then, for  any  $\de>0$ small enough one can choose $R$ so that
\beqn\label {Ga2}
P\Big[S^u[0,\tau_{D\setminus U_{ \de}}]\in \Ga\,\Big|\, S^u(\tau_{D})\in J \Big]  
> \frac{c}{2} P\Big[ W^u[0,\tau_{D\setminus U_{ \de}}]\in \Ga'\, \Big|\,  W ^{u}(\tau^W_D)\in J \Big].
\eeqn
If  $\de$ and $R$ are chosen independently of  $u$ and $D$ in (\ref{Ga}), then so are they in  (\ref{Ga2}).
(Here $S^u[k,n]= (S^u_j)_{k\leq j\leq n}$ and analogously for  $W^u[s,t]$.)
\end{Cor}
\v2\n
\pf\,  In view of the preceding corollary  it suffices to show
\beqn\label{eq_cor4.5}
P\Big[S^u[0,\tau_{D\setminus U_{ \de}}] \in \Ga \; \;\mbox{and}\;\; S^u(\tau_{D})\in J\,\Big] > \frac34 P\Big[S^u[0,\tau_{D\setminus U_{ \de}}] \in \Ga \; \;\mbox{and}\;\; S^u(\tau_{D\setminus U_\de})\in U_{\frac12 r, 2\de}\,\Big]. 
\eeqn
By Corollary  \ref{eqlem5.1} it follows that for all  $\de$ small enough, we can choose $R$ so that    if $w \in U_{\frac12 r, 2\de}$, then $P[ S^w(\tau_D)\in J] > 3/4$.  Hence  the right side of (\ref{eq_cor4.5}) is less than
 $P[S^u[0,\tau_{D\setminus U_{ \de}}] \in \Ga,  S^u(\tau_{D\setminus U_\de})\in U_{\frac12 r, 2\de} \; \;\mbox{and}\;\; S^w(\tau_D)\in J\,],$ which plainly entails  (\ref{eq_cor4.5}). \qed

\v2\v2\n
{\sc  4.3. Starting near the boundary (conditional case) }
\v2
 The main result of this subsection 
 (Proposition \ref{prop2}) concerns the walk conditioned to escape immediate absorption into the complement of $D$.
  It provides an  estimate of a conditional probability, given that the walk started at $v_0$ immediately enters into $D$ and hits $\partial U_{\La\de, \de}\cap D$ before leaving $D$.  Proposition \ref{prop2}, while  playing a crucial role in the next section,  is  not used in the succeeding subsections of the present section.

  Let $U$ be a domain of $\C$ whose boundary contains a vertex $v\in V$. The definition given in Section 2  of the first exit time $\tau_{U }$  for the walk starting from a boundary point of $U$ may be written as
  \begin{eqnarray}\label{tau_U}
\tau_{ U} &=& 1+ \tau_{ U}\circ \th_1 \qquad\quad \mbox{if} \quad S^v_1\in U, \nonumber\\
&=& 1 \qquad\qquad\qquad\quad\mbox{if} \quad S^v_1\notin U.
\end{eqnarray}
Here  $\th_n$ denotes the usual shift operator acting on random walk paths.  When there are more than two prime ends which are associated with  $v$,  the condition $S^v_1\notin U$ in (\ref{tau_U}) must be replaced by another one for the present purpose. Let  exactly $j$ prime ends, $v_k,\ldots,v_j$ say, correspond to $v$.
 Then  any continuous curve in  $U$ approaching $v$ may be considered to approach  one of these prime ends and not any other.  We say $u\in V$ is a neighbor site of $v_k$ in $U$ if $[v,u]$ is an edge of $G(U)$ and  the segment   $[u,v)$ approaches  $v_k$ ($k=1,\ldots,j$). Let ${\rm nbd}_U(v_k)$ denote the set of neighbor sites of $v_k$ in $U$. Then the latter condition in  (\ref{tau_U}) is replaced by
\beqn\label{tau_p}
\tau_U =1 \quad\quad\mbox{if} \quad S^{v_k}_{1} \notin {\rm nbd}_U(v_k)  
\eeqn
$ (k=1, \ldots, j)$ in order  to distinguish $v_k$ from the others.  %Although the arguments made in  the sequel can be adapted to either way, for definiteness as well as for simplicity of description
 In the sequel we adopt the latter definition  (\ref{tau_p})  with convention that
\v2
\quad {\it If we consider  $S^v$ with  $v\in V(\partial D)$,  the same letter $v$ is understood to designate \\
\quad\quad \;  a prime end of $\partial D$ that is associated with $v$.}
\v2
For typographical reason we often write $U(\La\de, \de)$ for $U_{\La\de, \de}$.

%Prop4.6
\begin{Prop}\label{prop2}\; For  any $\e>0$ there exists $\La>8$ such that  for   any $\de_0>0$ and any   $M>1$,  one can find $R>1$ such that if  $\de_0\leq  \de < 2/\La$,  $D \in \cD_{R,M}$, $v\in \overline{U}_{\de, \de }$ and $P[S^{v}(\tau_{ U(\La\de, \de)})\in  D]>0$, then
\beqn\label{J0}
 P[ \bar S^{ v}(\bar \tau_{ U(\La\de, \de)}) \in U_{\de} \, |\, S^{v}(\tau_{ U(\La\de, \de)})\in  D] < \e.
\eeqn
(The event under  $P$ says that the walk exits $U_{\La\de, \de}$ from either one of its two   narrow edges.; see (\ref{bar_tau}) for the notation $\bar \tau_U$.) 
\end{Prop}
\v2\n
\pf\;   Let  $U^+_\de$ denote  the upper   half  of the annulus  $U_{\de }$ which is  defined to be  the $\fa^{-1}$-image of   $\{\Im z>0\} \cap\fa( U_{\de })$. By symmetry it suffices to show
\beqn\label{J}
 P[\bar S^{ v}(\bar \tau_{ U(\La\de, \de)}) \in U_\de^+\, |\, S^{v}(\tau_{ U(\La\de, \de)})\in  D] < \e.
\eeqn

 The rest of the proof is broken into five steps. First we prove the proposition when $v=v_0$ in the  steps 1 through 4. The general case readily follows from this special case and is dealt with in the step 5. 
 
 \v2
 %%{\it Step 1}
{\it Step 1}. Put
$$T =\tau_{U(\La\de ,\de )}\quad\mbox{and}\quad \bar T =\bar\tau_{U(\La\de ,\de )},$$
and  let  $p(v)$ denote the conditional probability on the left-hand side of (\ref{J}):
\beqn\label{p00}
 p(v) = P[\bar S^{v}_{\bar T} \in U_\de^+\,|\,  S^{u}_{T} \in D], \qquad  v\in V(\overline{U}_{\La \de,\de}).
 \eeqn
 We claim that for $0< r < (\La-1)\de$, 
\beqn\label{y_v0}
 p(y)\geq p(v_0) \qquad \mbox{if \,  $P[S^y_T \in D]>0$ \,and \, $y\in V(\overline{U^+_{r,\de}}\cap \partial D)$},
 \eeqn
where $U^+_{r,\de} = U^+_{\de}\cap U_{r,\de}$. For the proof we use the fact that if  $X =(X_n)$ is
our walk $(S_n)$ killed at $T$ and conditioned on exiting through $\partial U_{\La\de ,\de }\setminus \partial D$ and if $h(v)= P[\bar S^v_T \in U^+_\de]$, then $X$  is the $h$-transform of our  walk $(S_n)$; in particular  the process $X$ is  Markovian and $p(v)$ is a harmonic function of it. Thus by Maximum principle applied  to $X$ there exists a path in $U_{\La\de ,\de }$ connecting  $v_0$ to $J:= U^+_\de \cap \partial U_{\La \de,\de}$ in $U_{\La \de,\de}$ such that $p(y) \geq p(v_0)$ for $y$ on the path.
If $\ga^{v_0}$ denotes such a path, then  the linear interpolation  $\bar S$ of the conditioned walk $X$  starting at $y\in \overline{U^+_{r,\de}}\cap \partial D$  must hit $\ga^{v_0}\cup J$ before exiting $U_{\La \de,\de}$ %since  $C^\circ\subset \overline{U_{\de\La,\de}\cap U^+_\de}$
 (recall $v_0$ is in a corner of $U^+_\de$), and the strong Markov property of $X$ concludes the claim (\ref{y_v0}).

 \v2
 %%{\it Step 2}
{\it Step 2}.  
 Let $U$ be a simply connected subdomain of $U_{\La\de ,\de }$ such that
 if 
 \beqn\label{C2}
 C=\partial U\cap \partial D,\quad \quad C^*= \partial U \cap  \partial U_\de\cap D,
 \eeqn
 then
$$\mbox{  \it  $C$ and $C^*$ are both  non-empty, $C$ is connected and}\quad C\subset \partial U^+_{\de}. $$
\v2\n
 Let  $v^*$ be a vertex in $C$.  In the next step $U$ and $v^*$ will be  specified more explicitly (by means of $D(x)$ defined in the preceding section), whereas in the present step  they may be rather arbitrary except for the restriction just mentioned.

  Put %$C^\circ =  \partial U\setminus \overline{\partial U \setminus C}= $
$$\Om =  \C\setminus \overline{\partial U \setminus C}.$$
  It is convenient to bring in
$$C^\circ = C\setminus \overline{\partial U \setminus C}.$$
Plainly $\partial  \Om =\overline{\partial U \setminus C} =\partial U \setminus C^\circ$ and $C^\circ$  
disconnects $U$ from $\Om \setminus \overline U$ in the graph $G(\Om)$.   Although the probability in question depends on our random walk $S^{v^*}$  restricted to  $D$, we are to extend it  to  $\Om$ through $C^\circ$. While how to choose the walk  outside  $U$ is at our disposal, we use the walk  $S^{v^*}$ itself  but with understanding that the extended walk  distinguishes  the same  vertex of  $ D$ according as  it is reached by  the walk  from a vertex of $C^\circ$   within  $\Om \setminus U$  or within $D$.  It  is 
% to define the \lq exit'  time $\tau_\Om$ in a  reasonable way in order 
 to distinguish two sides of $\partial \Om =\partial U \setminus C^\circ$ what we actually  need in the sequel, and   the  following notation  may  allow us to dispense with the   formal definition. Thus we write
 $$\lan C^*\ran_{out} =\{v\in V(D): [u,v] \in E \;\mbox{and}\; [u,v]\cap C^* \neq \emptyset \; \mbox{for some} \; u\in V(U\cup C^\circ)\}$$
(the subscript \lq$out$' reflects the fact  that  the directed segment $ [u,v]$ in the definition above is an \lq outward' boundary edge relative to $U$ although $U$ does not appear in the notation). By means of  this notation  the event that the walk exits from $\Om$ by crossing $C^*$ with an  edge directed outward from $U\cup C^\circ$ is expressed in the formula  $S^{v^*}_{\tau(\Om)} \in \lan C^*\ran_{out}$. We shall use  $\lan\partial U\ran_{out}$ in the analogous sense. 
 \v2
 %%{\it Step 3}
{\it Step 3}.  
Let $v^*\in V(C^\circ)$ and $L$ denote  the last exit time from   $\Om \setminus  U$ of the process $S^{v^*}_n,  n < \tau_\Om$: 
\[
 L =\left\{ \begin{array} {ll} \max\{0 \leq n < \tau_\Om: \;S^{v^*}_n \notin  U\} \quad &\mbox{if} \quad S^{v^*}_1 \notin U,\\
    0\qquad \qquad \qquad &  \mbox{if} \quad  S^{v^*}_1\in U,
 \end{array}\right.   
\]
so that $S^{v^*}_{L}\in C^\circ$ unless  $L+1= \tau_\Om$ since $C$ is a section cut from the boundary of the grid domain $D$.  
We are to compute 
$$q := P[S^{ v^*}_{L} \in C^\circ ,\; \bar S^{ v^*}_{\bar T}\circ \th_L \in U_\de^+,\,   S^{ v^*}_T\circ \th_L \in D], $$
the probability that the walk is found  in   $C^\circ$ at the epoch  $L$ and continued thereafter     and then  exits $U_{\La\de ,\de }$ through $U^+_\de \cap \partial U_{\La \de,\de}$ without landing on $\partial D$.
(Here    the shift operator $\th_L$ acts on $T$ as well as on $S^{v^*}$  as usual.)   
% our definition of $\tau_\Om$  is irrelevant to the  actual \lq  exit' from $\Om$ and the walk killed at time $\tau_\Om$ which is   in $C$ alive  at the present is never killed in the next step even if it steps  onto $\partial \Om$, showing  that    the event $S^{ v^*}_{L} \in C$ entails that the walk enters into $U$ in the  step next to $L$.   Then we 
Noting that if
$  n<\tau_\Om$, $S^{ v^*}_n\in C^\circ$  and  $\bar S^{ v^*}_{\bar T}\circ \th_n \in U_\de^+$, then  
 $L=n$,  we deduce that
\begin{eqnarray} \label{q1}
q&=&\sum _{n=0}^\infty\sum_{y\in C^\circ} P[ L =n <\tau_\Om, S^{ v^*}_n =y; \;\bar S^{ v^*}_{\bar T}\circ \th_n \in U_\de^+,\,  S^{ v^*}_T\circ \th_n \in D] \nonumber\\
&=&\sum _{n=0}^\infty\sum_{y\in C^\circ} P[ S^{ v^*}_n =y, n<\tau_\Om; \; \bar S^{ v^*}_{\bar T}\circ \th_n \in U_\de^+,\,  S^{ v^*}_T\circ \th_n \in D ]  \nonumber\\
&=&\sum _{n=0}^\infty\sum_{y\in C^\circ} P[S^{ v^*}_n =y, n<\tau_\Om] P[\bar S^{y}_{\bar T} \in U_\de^+, S^{y}_T \in D\,]   \nonumber\\
&=&\sum_{y\in C^\circ} G_{\Om}(v^*,y)P[S^{y}_{T} \in D]  P\Big[\bar S^{y}_{\bar T} \in U_\de^+\, \Big|\,  S^{y}_{T} \in D\Big],
\end{eqnarray}
where  $G_{\Om}(v^*,y) = \sum _{n=0}^\infty P[S^{ v^*}_n =y, n<\tau_\Om].$
Recalling the notation  $\lan C^*\ran_{out}$ (introduced   after $\Om$ is) we have 
$$P[S^y_T \in D]  \geq 
%\sum_{w\in \partial U\cap U_\de} P[S^y_{\tau(U)}=w]P[S^w_{T}\in D] + P[S^y_{\tau(U)}\in C^*].\geq
P[ S^y_{\tau(U)}\in \lan C^* \ran_{out}].$$
Hence
\beqn\label{4.19}
 q\geq \sum_{y\in C^\circ} G_{\Om}(v^*,y)P[S^{y}_{\tau(U)} \in \lan C^*\ran_{out}] p(y), 
\eeqn
 %where $T' $ is defined with $n\geq 0$ in place of $n\geq 1$ in the definition of $T$.  
where  $p(y)$ is the conditional probability defined in Step 1.

Noting that if $S^{ v^*}_{\tau(U)}\circ \th_L \in  \lan C^* \ran_{out}$  and $S^{ v^*}_{L} \in C^\circ$,  then $S^{ v^*}_{\tau(\Om)} \in \lan C^* \ran_{out}$ and vice versa,  we have
\beq
 \sum_{y\in C^\circ} G_{\Om}(v^*,y) P[S^y_{\tau(U)} \in  \lan C^* \ran_{out}] 
 &=& \sum _{n=0}^\infty\sum_{y\in C^\circ} P[ n<\tau_\Om, S^{ v^*}_n =y;\; S^{ v^*}_{\tau(U)}\circ \th_n \in  \lan C^* \ran_{out}]\\
 &=& \sum _{n=0}^\infty P[ L=n, S^{ v^*}_n \in C^\circ;\; S^{ v^*}_{\tau(U)}\circ \th_n \in  \lan C^* \ran_{out}]\\
% &=& \sum _{n=0}^\infty P[L= n;\; S^{ v^*}_{\tau(U)}\circ \th_n \in  \lan C^* \ran_{out}]\\
&=&  P[S^{ v^*}_{\tau(\Om)} \in \lan C^* \ran_{out}].
\eeq
Thus in view of (\ref{y_v0}) and (\ref{4.19})
 \begin{eqnarray}\label{5.3}
 q \geq   p(v_0) P[S^{v^*}_{\tau(\Om)}\in \lan C^* \ran_{out}].
 \end{eqnarray}

%%%%%%
%{\it Step 4}
\v2
{\it Step 4}.  
In this step we verify   the inequality of the proposition for $v=v_0$. First we  claim  that one can choose  the pair  of $U$ and $v^*$  so that  for some universal constant $c^\circ>0$
 \beqn\label{2.4}
 P[S^{v^*}_{\tau(\Om)}\in  \lan C^* \ran_{out}] \geq  c^\circ \quad \mbox{ if}\; \La\geq m
 \eeqn
 and that $U$  satisfies   (along with those stated at (\ref{C2}))
 \beqn\label{2.41}
 U \subset  U^+_{2m\de,\de}.
  \eeqn
Here $m$ is  the universal constant specified in Lemma  \ref{LSW}.  To this end we take  up $U$  as given by Lemma \ref{LSW2} but with (\ref{2.41}) in place of the inclusion $U\subset U_{m\de,\de}$ (note that   $\fa(U^+_{2m\de,\de})$ is simply a rotation of $\fa(U_{m\de,\de})$ and the special choice of $x$ as $\fa^{-1}(1-\de)$    in the proof carries no significance except for this inclusion).  
As $v^*$ we take a  vertex in $\partial D$ closest to $x^*\in \partial U$ (described in Lemma \ref{LSW2}). Then an easy application of Corollary \ref{YY4.2}  verifies (\ref{2.4}) with $c^\circ=c^*/2$, provided $R$ is large enough and  $D\in \D_{R,M}$.   Plainly we have (\ref{2.41}). Thus the claim has been  proved.

From the last expression of $q$ in (\ref{q1}) and by using the last exit decomposition as in Step 3 (but in reverse direction) we have  
\beq
q&\leq & \sum_{y\in C^\circ} G_\Om(v^*,y) P[\bar S_{\bar T}^y \in U_\de^+ ]\\
&\leq&\sum_{y\in C^\circ} G_\Om(v^*,y)\sum _{u\in \lan \partial U\ran_{out}\cap U_\de}  P[S^{y}_{\tau(U)}  =u\,] P[\bar S_{\bar T}^u \in U_\de^+]\\
&\leq& P[ S^{v^*}_{\tau(\Om)}\in \lan \partial U\ran_{out}\cap U_\de] \sup_{u \in \lan \partial U\ran_{out}\cap U_\de} P[\bar S_{\bar T}^u \in U_\de^+] \\
&\leq&  \sup_{u\in \lan \partial U\ran_{out} \cap U_\de} P[\bar S_{\bar T}^u \in U_\de^+],
\eeq
and by  repeated applications of  Corollary \ref{cor2.2} we conclude  $q\leq c_1^{\lfloor\La - 2m\rfloor}$ with a universal constant  $c_1<1$ ($\lfloor a \rfloor$ denotes the largest integer that does not exceed a real number $a$). From (\ref{5.3}) and (\ref{2.4})   it therefore follows   that $p(v_0) \leq c_1^{[\La-2m]}/c^\circ$,  showing  (\ref{J})  for $v=v_0$. % with the same $\La$ as in Step 3.

%%%%%%
%{\it Step 5}
\v2
{\it Step 5}. Consider the general case $v\in \overline{U}_{\de,\de}$. We must prove that $p(v)$, the probability defined  by (\ref{p00}), can be made arbitrarily small by taking $\La$ large enough.   If $v\in \partial D$,  the same proof as above  apply since $1\in \partial \D$ may be replaced by any point of  $\partial \D$. 
 Let $v\in D \cap \overline{U}_{\de,\de}$.  We suppose  $v\in U^-_\de =\fa^{-1}(\{\Im z<0\})\cap U_\de$ for simplifying the description.      Let $\ga$  be a path from $v$ to $J$  in $U_{\La\de,\de}$ on which $p(v)\leq p(u)$, where $J=U^+_\de \cap \partial U_{\La \de,\de}$ as before.       Let $w$ be a vertex on $\partial D $ closest  to
$\fa^{-1}(e^{i \La\de/2})$ and $A$ the event that the walk $S^{w}$ exits $U^+_{\La\de,\de}:= U_{\La\de,\de}\cap U^+_\de$ without hitting $\ga$. Then,  noting that  the event $S^w_T\in D$ (with $T= \tau_{U_{\La\de,\de}}$ as before) entails  
$S^{w}(\tau_{ U^+_{\La\de,\de}})\in  D$ and writing $B$ for the latter event, we infer that 
$$ P[A \, |\, S^{w}_T\in  D] 
=  \frac{P[A, \, S^{w}_T\in  D\, |\, B] } {P[ S^{w}_T\in  D\, |\, B]}
\leq  \frac{P[A\, |\, B] } {P[ S_T^{w}\in  D\, |\, B]}.
$$
  In view of what is noted at the beginning of this step the last ratio as well as $p(w)$ may be made arbitrarily small by taking $\La$ (and $R$) large enough.  On the other hand on using   the strong Markov property of the walk conditioned  on $S^w_{T}\in D$ (as in Step 1) we deduce the inequality
  $$p(w) \geq P[ \ga \;\mbox{is hit before exiting}\; U^+_{\La\de,\de}, \bar S_{\bar T}^w \in U^+_\de\,|\,  S^{w}_T\in  D]  \geq (1-P[A \, |\, S^{w}_T\in  D] )p(v).$$
   Hence   $p(v)$ may be made arbitrarily small.  
  The proof of Proposition \ref{prop2} is complete. \qed

\v2

Taking $\e=1/2$ in Proposition \ref{prop2} we plainly obtain the following
%Cor4.7
\begin{Cor}\label{cor2}  For   any $\de_0>0$ and  $M>1$,  one can find $R>1$ such that if $\de_0\leq \de <1/4$,  $(\hat o,D) \in \D_{R,M}$,  $v\in V(\overline{U}_{\de})$,  $P[S^{v}(\tau_{ U(\de, \de)})\in  D]>0$ and $|\Im \fa_D(v)|<\de$, then for some universal integer $m^*\geq 4$,
\beqn\label{cor5.2}
 P[ \bar S^{ v}( \bar \tau_{ U(m^*\de, \de)}) \in  \partial U_{\de} \cap D \, |\, S^{v}(\tau_{ U(m^*\de, \de)})\in  D] > 1/2.
\eeqn
\end{Cor}

\v2\n
{\sc 4.4. \; Hitting distribution of $\partial D$}
\v2
 This and the next subsections, in which we do not use Proposition  \ref{prop2},   primarily concern   the hitting distribution
  $$H_D(u,b) :=P[S^{u}_{\tau(D)} =b],  \quadd u\in V(D),  \; b\in \partial D,$$
and provide some estimates of it.    In later applications we need extend it to $b\in V(D)$ by
$$H_D(u,b) :=P[S^{u} \, \mbox{ visits $b$ before exiting $D$}\,],  \quadd u\in V(D)\setminus \{b\}.$$
    If $u\in V(\partial D), u\neq b$ and $\hat u$ is a prime end that is associated with $u$, we set 
    $$H_D(\hat u,b) = \sum _{v\in {\rm nbd}_D(\hat u)} P[S^{\hat u}_1=v] H_D(y,b)$$
     to be consistent to the definition of $\tau_U$ in (\ref{tau_U}).  In the sequel, however,  we write simply $v$ for $\hat v$ according to the convention advanced right
     before Proposition \ref{prop2}. 

The most results presented below are essentially the same as what are  found  in \cite{YY}.
We give  proofs to some of them, which are simpler than those in  \cite{YY} mainly owing to Corollary \ref{YY4.2},   although the idea of the proofs are the same as in \cite{YY}.  In the  proofs we shall often drop   the subscript $D$ from  $H_D$ as well as from $\fa_D$. 

Given $r>0$ and $b \in V(\overline{U}_{r /4})$, put 
$I= \{e^{i\th}\fa_D(b): |\th|\leq r\}$,  $J= \fa_D^{-1}(I)$. 
 The next result is essentially  the same as Lemma 5.8 of \cite{YY}.  
%Lem4.8
\begin{Lem}\label{lem4}\; For  any    $0<r<\frac1{10}$ and  $M>1$   there exists $R$  such that if  $D\in \D_{R.M}$ and  $b \in V(\overline{U}_{r /4})$,  then for all $v\in V(D\setminus U_r)$ and  $w\in V(D\setminus B_{4r}(b))$ with $P[ S^{w}(\tau_D)\in  J]>0$, 
$$P[ S^{w}(\tau_{D})=b \, |\, S^{w}(\tau_{D}) \in J]\leq c_2P[ S^{v}(\tau_{D})=b \, |\, S^{v}(\tau_{D}) \in J]$$
for some universal  constant  $c_2$.
\end{Lem}
\n
\pf\; For simplicity we  suppose  $b\in \partial D$, the arguments below being readily adapted to the case $b\in V(\overline{U}_{r /4})$.   Put $q_J(y)= P[ S^{y}(\tau_{D})=b \, |\, S^{y}(\tau_{D}) \in J]$, $y\in V(D)$.  It suffices to prove that if 
% (see (\ref{angle}) for the notation $\lan K\ran$)  
    $v\in V(D\setminus U_r)$, then  for some universal  constant  $c^\circ >0$, 
\beqn\label{eq5.2}
q_J(y)\leq c^\circ q_J(v) \quad \mbox{for}\quad y \in  V(B_{4r}(b)\setminus \overline {B_{3r}(b)}\,).
\eeqn
For   $q_J(w)$ ($w\notin B_{4r}(b)$)  is a convex combination of $ q_J(y)$,  $y \in  V(B_{4r}(b)\setminus \overline{ B_{3r}(b)}\,)$, provided  that $R$ is so large that  no  edge in $G(D)$  joins $\overline B_{3r}(b))$ and $D\setminus B_{4r}(b))$.   
 
 Since $q_J$ is harmonic for the walk conditioned to exit $D$ through $J$, there exists a path $\ga^y$ of $G(D)$ joining $y$ and $b$  such that $q_J(y)\leq q_J(z)$ for all $z\in V(\ga^y)$.  Now let $v\in V(D\setminus U_r)$ and consider the event, denoted by  $B$,  that   the path $(S^v_n)_{0\leq  n <\tau_D}$  enters into   $\overline{B_{2r}(b)}$ avoiding  $U_{r/2}$ as the landing place: 
 $$
 B = \{ S^v( \tau_{D\setminus B_{2r}(b)}) \in D\setminus U_{r/2} \}.
 $$
 %(which always occurs if   $v\in \overline {B_{2r}(b)}$).  
 We   apply  Corollary \ref{ncor4.5} and Corollary \ref{YY4.2}  to see that 
  given  the event  $B$ as well as   the event $S^{v}(\tau_{D}) \in J$  occurs,
   the conditional  probability that the  walk  $S^v$ crosses $\ga^y$ before exiting $D$  is bounded below by a universal  positive constant. 
   (To  this end one may  first verify   that  under this conditioning the conditional probability that for any $\eta>0$ small enough, the walk starting at a vertex in $ B_{(2+\eta)r}(b) \setminus U_{r/2}$ exits $B_{3r}(b)$  through the upper  half of  $J$ is bounded
   below by a universal constant, and by symmetry the same is true for the exiting through the lower half.)   
On the other hand,   applying  Corollaries \ref{ncor4.5} and \ref{YY4.2} again we infer that 
 $P[B\, |\, S^{v}(\tau_{D}) \in J] \geq c_1$  with a universal constant $c_1>0$.
   Hence,  the  conditional  probability of  the  walk  $S^v$ crossing $\ga^y$ before exiting $D$  given $S^{v}(\tau_{D}) \in J$  is bounded below by a positive  universal constant, $c_*$ say, from which, on using the strong Markov property of the conditioned walk,   we infer 
 $$q_J(v) \geq P[S^v(\tau_D)=b, S^v[0,\tau_D] \cap \ga^y \neq \emptyset\, |\, S^{v}(\tau_{D}) \in J]  \geq c_* q_J(y),$$
showing   (\ref{eq5.2}) with $c^\circ =1/c_*$.   This finishes the proof of  Lemma \ref{lem4}. \qed

\v2\n

%cor4.9
\begin{Cor}\label{prop5.2}    For  any    $0<\de<\frac1{10}$ and  $M>1$   there exists $R=R(\de)$  such that if  $D\in \D_{R.M}$ and    $b \in V(\overline{U}_{\de /4})$,  then for all $v\in V(D\setminus U_\de)$ and  $w\in V(U_{2\de}\setminus B_{1/2}(b))$, 
\beqn\label{ieq_p5.2}
 H_D(w,b)
\leq c_3 H_D(v,b)
\eeqn  
with a universal constant $c_3>0$.
\end{Cor}

\v2\n
\pf \;  Let  $I$ and  $J$ be as in Lemma \ref{lem4} with $r=\de$  and $R$ be chosen large enough. 
 Then according to Corollary \ref{cor5.3}, uniformly for  $D\in \cD_{R,M}$ and  $v\in D\setminus U_\de$,
\beqn\label{H_h1}
H(v;J) = P[S^v(\tau_D) \in J] \geq \frac12  P[ W^{\fa(v)}(\tau_{\D})\in I] \geq h(\de),
\eeqn
with  $h(\de)= c\de^2$ for some $c>0$   (one may take $c= 1/8\pi$). 
According to Corollary \ref{lem2.2} we can find a number $\eta=\eta(\de)$ such that $ P[S^w(\tau_{B_\de(w)})\in D] \leq h(\de)$ for $w\in U_{\eta\de}$.  Hence if $ w\in U_{\eta\de} \setminus B_{3\de}(b)$, then 
 \beqn\label{H_h2}
 H(w;J) \leq  h(\de). 
 \eeqn
On the other hand  for $w\in (U_{2\de}\setminus U_{\eta\de})\setminus B_{1/2}(b)$, we see 
$P[ W^{\fa(w)}(\tau_{\D})\in I] \leq 2c^\circ h(\de)$ where $c^\circ$ is a universal constant. Thus 
the bound  (\ref{H_h2}) with $h(\de)$ replaced by $c^\circ h(\de)$ is  valid also for such $w$.    Combined  with  (\ref{H_h1}) we then conclude that if $w\in U_{2\de}\setminus B_{1/2}(b)$ and $v\in D\setminus U_\de$,  
$$\frac{H(w,b)}{H(v,b)} =\frac{P[S^w(\tau_D)=b\,|\, S^w(\tau_D)\in J]}{P[S^v(\tau_D)=b\,|\, S^v(\tau_D)\in J]} \frac{H(w; J)}{H(v,J)}\leq  c_2 c^\circ,$$
where $c_2$ is the constant in Lemma \ref{lem4}. The proof of Corollary \ref{prop5.2} is complete.

\v2 
{\sc Remark 3.}  The estimate of   Corollary  \ref{prop5.2} holds also for any pair $w\in U_\de$ and $v\in D\setminus U_\de$ such that $\fa_D$ sends both  $w$ and $v$ into $\{z\in \D: k\de<|z-\fa(b)|\leq (k+1)\de\}$ for some $k=3,4,\ldots$.

\v2\v2\n
{\sc  4.5.  Poisson kernel approximation }
\v2
Let $\la$ be the usual Poisson kernel for $\D$:
 $$\la(z,w) = \frac{1-|z|^2}{|z-w|^2} \quad  (z\in \D,  w\in \partial \D).$$
The following proposition is proved in \cite{LSW} when $G$ is the square lattice and in \cite{YY} in the same   setting as ours. %We give a proof which is  based on  the same idea of  but much simpler than that of  \cite{YY}.

 %Prop4.10
  \begin{Prop}\label{prop5.3}  For any $\e>0$ and $M>1$  there exists $R_0=R_0(\e,M)$ and  $\eta>0$ such that if  $D\in \cD_{R,M}$, $b\in V(\partial D)$ and   $a\in V(D\setminus U_\e)$  ( i.e., $a\in V(D)$ and $| \fa_D(a)|<1-\e$),   then  for $b' \in \overline{B_\eta(\fa_D(b))}$ with $H_D(\hat o, b')>0$, 
 $$\bigg| \frac{H_D(a,b')}{H_D(\hat o, b')} - \la(\fa_D(a),\fa_D(b)) \bigg|\leq { \e}. $$
 \end{Prop}

For our application in the next section we work with the upper half plane $\H =\{z; \Im z>0\}$ and it is convenient to translate the formula of Proposition \ref{prop5.3} as in the  next corollary.
%Cor4.11
\begin{Cor}\label{cor5.30}
For any  constants $\epsilon >0$ and $M>1$, and any compact set  $K$ of $\H$,  there exists  $R>0$ such that
if $D \in \cD_{R,M}$,   $\psi  : D\rightarrow \H$ is a conformal map with $\psi (\hat o)=i$ and  $b \in V(\partial D)$ with $H_D(\hat o,b)>0$, then for all $y$ and  $w$ from $V(D)$ 
 such that $\psi(y), \psi(w) \in K $, 
 \[\left| \frac{H_D(w,b)}{H_D(y,b)}-\frac{\, |\psi(y)-\psi(b)|^2\Im \psi(w)\,}{\,|\psi(w)-\psi(b)|^2 \Im \psi(y)\, }\right|<\epsilon.\]
\end{Cor}

\section{ Convergence of LERW to chordal SLE$_2$}

This section concerns the convergence  to a chordal SLE$_2$ of 
  the loop erasure of the random walk on the planar graph $G$ started at a boundary vertex of a grid domain $D$ and conditioned to exit at another boundary   vertex.  After giving a brief exposition  of the chordal Loewner chain together with  a few preliminary lemmas  in Subsections 1 and 2 we state our result on the convergence of LERW  (Theorem \ref{main1}) and advance an abridged proof of it  in Section 5.3. 
  \v2
  
{\sc 5.1. Chordal Loewner Chain for a simple curve  in $\H$}
\v2

A chordal Loewner chain is the solution of a type of Loewner equation that describes the evolution of a continuum  growing 
from the  boundary to the  boundary of a simply connected domain of $\C$. For our present purpose we have only to consider   the case when the continuum is a  simple curve. 
In this subsection we consider the special case when the domain is $\H:=\{z \in \C : \Im\,  z>0 \}$,  the upper half plane
and the curve grows from the origin to the infinity in $\H$, general case will be  considered in the next subsection.

Suppose that $\gamma :[0,\infty) \rightarrow \overline {\H}$ is a simple curve 
with $\gamma (0)\in \R,\gamma (0,\infty) \subset \H$.
Then, for each  $t\geq 0$, there exists a unique conformal map 
$g_t:\H \setminus \gamma(0,t] \rightarrow \H $ 
satisfying $g_t(z)-z \rightarrow 0$ as $z \rightarrow \infty $. It is noted that  $g_t$ can be continuously extended 
to the (two sided) boundary of  $\H\setminus \gamma (0,t]$ along $\gamma (0,t]$.  For each $t$, there exists the limit
\[ {\rm hcap}(\ga[0,t]):=\lim _{z \rightarrow \infty} z(g_t(z)-z), \]
called the  half-plane capacity of $\ga[0,t]$; ${\rm hcap}(\ga[0,t])$ is real, and increasing and continuous in $t$.
If $\gamma$ is parametrized by $\frac12$ times  the  half-plane capacity  (so that  $ {\rm hcap}(\ga[0,t])=2t$), then according to Loewner's theorem
$g_t$ satisfies his differential equation
\begin{equation} \label{cle}
\frac{\partial }{\partial t}g_t(z)  =\frac{2}{g_t(z)-U(t)},
\quad g_0(z)=z, 
\end{equation}
where $U(t)=g_t(\gamma(t))$ and $U(t)$ is a $\R$-valued continuous function (see \cite{L}).
The equation (\ref{cle}) is called the {\it chordal Loewner equation} and $U(t)$ the {\it driving function}. The family  $g_t,  t\geq 0$ is called the  {\it chordal Loewner chain} generated by  a curve $\ga$ and/or driven by a function $U(t)$.

Conversely, given a continuous function  $U(\cdot) :[0,\infty)\rightarrow \R$, one can    solves   the ordinary differential equation (\ref{cle})  for each  $z \in \H$ to obtain the solution $g_t(z)$ up to the time
$T_z:=\sup\{t>0:|g_t(z)-U(t)|>0\}$. A function $\ga$ may be defined by $\ga(t)= \lim_{z\to U(t), z\in \H} g_y^{-1}(z)$, provided the limit exists. If $U(t)$ has a sufficient regularity, this gives a simple curve $\ga[0,t] =\{z \in \H:T_z\leq t\}$ and then for $t>0$, $g_t(z)$ is a conformal map from $\H \setminus \ga[0,t] $ onto $\H$.
%The family $(g_t)_{t \geq 0}$ describes the evolution of hulls $(K_t)_{t \geq 0}$ corresponding to $U(\cdot)$ and growing from the boundary to $\infty$.
%Therefore, we have a one-to-one correspondence between $U(\cdot )$ and $(K_t)_{t \geq 0}$.
If $U(\cdot)$ is the driving function of a simple curve $\gamma $ in particular, we can recover $\gamma $ from $U(\cdot)$. If this is the case  the curve $\ga$ may be said to be driven by $U$ (via the Loewner chain).

It is known that  if we take   a linear Brownian path $\sqrt{\k}\, W(t)$   as $U(t)$, a simple (random) curve $\gamma $ is driven by $U$  for $0<\k\leq 4$ with probability one; in the case $\k > 4$, the procedure above still produces a curve which is, however, no longer a simple curve. In either case the random curve driven by $\sqrt{\k}\,W$ is called a {\it  chordal $SLE_\k$} curve in $\H$.

%If $U(t)$ is sufficiently nice, then $K_t$ is generated by a curve $\gamma$ 
%with $\gamma (0) \in \R, \lim_{t\to\infty} \gamma (t) =\infty$
%(i.e., $\H\setminus K_t$ is the unbounded component of $\H\setminus \gamma(0,t]$).
%However, there exists a continuous function $U(t)$ such that
%$K_t$ can not be generated by a curve.
%There is known a sufficient condition for $U(t)$ to drive a curve as given by  
%\begin{pp}  \label{gen_sim}
%Suppose  for some $r < \sqrt{2}$ and all $s < t$,
%\[|U(t)-U(s)| \leq r\sqrt{t-s} .\]
%Then $K_t$ is generated by a simple curve.
%\end{pp}

%In summary, a simple curve $\gamma$ brings out a Loewner chain, whereby it determines the driving function $U(t)$, and conversely a continuous function $U(t)$ with appropriate regularity generates a simple  curve  through the Loewner chain driven by $U(t)$.

%%%%%
\v2
In the rest of this subsection  let  $U(t)$ be  a real continuous function of $t\geq 0$. We  shall need the following  lemma, which   is   Lemma 2.1  in \cite{LSW} but restricted to the case when  a simple curve is driven by $U$. 
%Lem5.1
\begin{Lem}  \label{k-est2}   Suppose that a simple curve $\ga(t)$ is   driven by  $U(t)$.  Put  
\beqn\label{kt}
k(t) =\sqrt{t}  + \max_{0\leq s\leq t} |U(s)-U(0)|.
\eeqn
Then, for any $t>0$,
$c^{-1}k(t) \leq {\rm diam} (\ga[0,t]) \leq c k(t)$ for   a universal constant $c>0$.
\end{Lem}

\v2
  From the Loewner equation (\ref{cle})  we have
\begin{equation} \label{ineq_Leq}
|g_t (z)-z| < t \cdot \sup _{0  < s \leq t} \frac{2}{|g_s  (z)-U(s)|},
\end{equation}
which  the upper bound  (easier half)  in Lemma \ref{k-est2} is deduced from. While the next lemma   improves  this upper  bound,  our application of it concerns   its another aspect.
%Lem6.2
\begin{Lem}\label{k-est1}   Put for $t>0$ and $z\in \H$
  $$\la_t = \min_{0\leq s\leq t}|z- U(s)| \;\; \mbox{ and}   \;\; \mu_t = \min_{0\leq s\leq t}|g_s(z)- U(s)|.$$ 
Then $\mu_t> \frac12\Big(\la_t+\sqrt{\la_t^2-8t}\Big)$ as long as $\la_t \geq \sqrt{8t}$.
%$ \min_{0\leq s\leq t} |z-U(s)| \geq \sqrt{8t}$  for $0\leq t\leq t_0$.  Then, $\min_{0\leq s\leq t_0}|g_s(z)-U(s)| >  \frac12\Big(\la_*+\sqrt{\la_*^2-8t_0}\Big)$, where  $\la_* = \min_{0\leq t\leq t_0} |z-U(t)|$. 
\end{Lem}
\v2\n
\pf\;  The proof rests only on (\ref{ineq_Leq}).  We deduce from it that $\la_t -\mu_t   <  2t/\mu_t$. Indeed this inequality certainly  holds true for  $t_*$ at which  the minimum  $\mu_t$ is attained so that $\mu_t=\mu_{t_*} = |g_{t_*}(z)-U(t_*)|$; hence it does also for $t$.
We may rewrite  it as  $\mu_t^2 - \la_t\mu_t + 2t >0$.  If   $\la_t\geq \sqrt{8t}$, then,  putting  
$\xi_\pm(t)= \frac12 \Big( \la_t \pm \sqrt{\la_t^2 - 8t}\Big),$  we have either  $\mu_t <\xi_-(t)$ or $\mu_t> \xi_+(t)$.  
%These  alternatives remain true if $\xi_+(t)$ and  $\xi_-(t)$ are replaced by
%$$\tilde\xi_+(t) := \frac12 \inf_{0\leq s\leq t}\Big( u_s + \sqrt{ u_s^2-8s}\Bigg)~~\mbox{and}~~ \tilde\xi_-(t) := \frac12 \sup_{0\leq s\leq t}\Big( u_s - \sqrt{  u_s^2-8s}\Bigg),$$
%respectively.
Clearly $\xi_-(t) \leq  \xi_+(t)$ and   $\mu_t-\xi_+(t)$ is continuous and positive for $t$ small enough, hence we must have $\xi_+(t) <\mu_t$ for all  $t\leq t_0$, where  $t_0=t_0(z)$ is determined by $\la_{t_0}= \sqrt{8t_0}$.    \qed
\v2

The expression  for  the logarithmic derivative of the  imaginary part of $g_t$ derived from   (\ref{cle}) we deduce that  for $z\in \H$, $t>0$,
\begin{equation} \label{ineq_Leq2}
 \frac{\Im\,  g_t(z)}{\Im\, z} 
  = \exp \left(- \int_0^t \frac{2}{|g_s(z)-U(s)|^2} ds\right).
\end{equation}
%Since $\min_{0\leq s\leq t}|z-U(s)|>|z-U(0)| - k(t) +\sqrt{8t}$,
 Lemma \ref{k-est1} combined with  (\ref{ineq_Leq}) and (\ref{ineq_Leq2}) entails the following corollary.

%%Cor6.3
\v2\begin{Cor}\label{cor_k-est}   %Let $k(t)$ be defined by  (\ref{kt}) and put $\la_t = \min_{0\leq s\leq t} |z-U(s)|$. If $|z-U(0)|\geq k(t)$, then 
If $\la_t\geq \sqrt{8t}$, then
$$  |g_t(z)-z|<\frac{4t}{\la_t+\sqrt{\la_t^2-8t}} \quad \mbox{and} \quad 1 \geq \frac{\Im g_t(z)}{\Im z} > \exp\bigg(-\frac{8t}{(\la_t+\sqrt{\la_t^2-8t}\,)^2 }\bigg).$$
\end{Cor}

% ,  since the assumption of Lemma \ref{k-est1}  is implied by the inequality 
% $|z|\geq \sqrt{8t_0} + \sup_{0<s<t_0}|U(s)|$. it shows that $\ga[0,t]$ is contained in the set
%$\{z\in \H: |z-U(s)|\leq \sqrt{8s}$ for some  $s\in [0,t] \}$, and in particular

\v2
\v2\n
%%%\subsection{6.2. Chordal Loewner chain in simply connected domains}
{ \sc 5.2. \;Chordal Loewner chains in simply connected domains}
\v2

We adapt  the formulation of \cite{SS}.  Let $D$ be a simply connected  domain  with  two distinct boundary points $v_0, v_e \in \partial D$   (to be precise these should be prime ends).  
Let $\psi : D \rightarrow \H$ be a conformal map with $\psi  (v_e)=0, \psi  (v_0)=\infty$.
Although $\psi$ is not unique, any other such map  can be written as $y\psi $ for some $y>0$, so that $v_0, v_e$ being given, the map $\psi$ is in one-to-one correspondence to  a point $\hat o$ on the curve 
$(\psi^{-1}(iy))_{0<y<\infty}$ 
(running  from $v_e$ to $v_0$ in $D$) via  the equation  $\psi(\hat o)=i$. We shall take this $\hat o$ as a reference point attached to $D$. 

For a simple curve  $\gamma :(0,T) \rightarrow D$  connecting  $v_e$ and $v_0$ so that  $\ga(0+)=v_e$ and $\ga(T-)=v_0$,
let $g_t$ be the Loewner chain generated by the curve  $\psi \circ \gamma : (0,T) \rightarrow \H$. Put
$D(t):= D\setminus \ga(0, t]$ and define  $\psi_t: D(t)\mapsto \H$ by  
\beqn \label{psi_t}
\psi_t = g_t\circ\psi|_{D(t)}, \quad t\in [0,\infty),
\eeqn
where $\psi |_{D(t)}$  designates the restriction of $\psi$ to $D(t)$. Now we  reparametrize  the curve $\psi\circ\gamma$ in $\H$ by half plane capacity  so that $2t= {\rm hcap}(\psi \circ \gamma [0,t] ), 0\leq t <\infty$.     The driving function $U(t)$ of the chain $g_t$ is then given by
\[U(t) = \psi_t(\gamma (t)).\]
The family of conformal maps  $\psi_t, t\geq 0$ may also be called a {\it chordal Loewner chain} (in $D$) with driving function $U(t)$.
For each $s>0$,  the curve $\ga^{(s)}(t) := \ga(s+t)$ connects $\ga(s)$ and $v_0$ in $D(s)$,  and $\psi_s$  conformally maps $D(s)$ onto  $\H$  with $\psi_s(\ga(s)) = U(s)$, 
$\psi_s(v_e)=\infty$.
On putting
\begin{equation*} \label{U1}
 g^{(s)}_t  = g_{s+t} \circ g_s^{-1} \quad \text{and} \quad \psi_t^{(s)} = \psi_{s+t}|_{D(s)},
 \end{equation*}
 substitution into  $U(s+t)=  \psi_{s+t}(\gamma (s+t))$ yields 
\begin{equation} \label{U2}
 U(s+t) = \psi_t^{(s)}(\ga^{(s)}(t)).
\end{equation}
It  follows that   $\psi_t^{(s)} = g_t^{(s)}\circ \psi_s$,  $g_t^{(s)} $ (and $\psi _t^{(s)}$) is the Loewner chain generated by  the curve  $\gamma^{(s)}$
and  $U^{(s)}(t) := U(s+t)$ is the driving function of the chain  $\psi _t^{(s)}$ in $D(s)$. 
    
Define $\hat o(t)\in D(t)$ by
\beqn\label{p(t)}
 \psi_t(\hat o(t)) = U(t)+i. 
 \eeqn
Then  $\hat o=\hat o(0)$,  and $\hat o(t)$ will serve as  an appropriate reference point of $D(t)$ for our purpose.

The following lemma, proved in \cite{Su},  plays a significant role to show Theorem \ref{main1} below. 
%Lem6.4
\begin{Lem} (\cite[Lemma 2.4]{Su}) \label{refpt}  Let $\hat o(t)$ and  $\rho_D(t)$  be as above.  
Given  $T>1$ and $\epsilon >0$,  put $\tilde T:= \sup\{t\in [0,T]: |U(t)|<1/\epsilon\}$. 
There then exists a constant $c(T,\epsilon)>0$, 
which    does not depend on $(\gamma(t), t\geq0)$ nor on $D$, such that $\mbox{in-rad}_{\hat o(t)} D(t)\geq c(T,\epsilon) \,\rho_D$  for  $t<\tilde T$.
\end{Lem}  
\n\pf \;   In \cite{Su} the constant $c(T,\e)$ is allowed to  depend on $(D, \ga(0))$, which is not important therein since  $D$ is 
fixed in the setting of its main theorem.  
In the  first half  of  the proof of Lemma 2.4 of \cite{Su}, its substantial  part,  it is verified    that 
 \beqn\label{claim1}
 |\phi(\hat o(t)) - \phi(\ga(t'))| \geq 2^{-1} e^{-4\tilde T}  \quad \mbox{if}\;\;  t'\leq t  <\tilde T.
 \eeqn   
We write down the  other half (with slight modification of wording) to ensure that $c(T,\e)$ can be taken independently of  $(D, \ga(0))$.
In \cite{SS} (the proof of Corollary 4.3) it is proved that 
$$\mbox{in-rad}_{\hat o(t)}(D) \geq c_0(T,\e)\mbox{in-rad}_{\hat o(0)}(D) \quad \mbox{and}\quad {\rm dist}(\phi(\hat o(t)),\R) \leq M(T,\e).$$ 
According to the first inequality     it  suffices for the proof of the lemma to show that
$$ |\hat o(t) -\ga(t')| \geq c_1(T,\e){\rm dist}(\hat o(t), \partial D) \qquad t' \leq t <\tilde T.$$
If $|\hat o(t) -\ga(t')| <2^{-1}{\rm dist}(\hat o(t),\partial D)$, we may take $1/2$ for $c_1(T,\e)$;  otherwise,  applying (\ref{claim1})  and  the distortion theorem
in turn yields
$$2^{-1}e^{-4\tilde T}\leq  |\phi(\hat o(t))-\phi(\ga(t')) | \leq 16|\hat o(t)-\ga(t')| \frac{{\rm dist}(\phi(\hat o(t)),\R)}{{\rm dist}(\hat o(t), \partial D)}. $$
 Hence
$|\hat o(t)-\ga(t')| \geq [e^{-4T}/32 M(T,\e)] {\rm dist}(\hat o(t), \partial D)$
as desired. \qed

%{\sc 5.3. ~LERW and its convergence to chordal  SLE}
\v2\v2\n
{\sc 5.3. \; Convergence of driving function}
\v2
Let $G=(V,E)$ be a planar irreducible graph whose edges are directed and weighted as   described in Section 2.
For any finite sequence $\om = (\om_0,\ldots, \om_N)$ in $V$, the {\it loop erasure} of $\om$, denoted by LE$(\om)$,  is the sequence obtained by erasing  the loops in it in chronological order.  To be precise LE$(\om)= (\hat\om_0,\ldots, \hat\om_L)$ is defined as follows: on  putting  $\hat\om_0 = \om_0$, for $k=0,1,\ldots$,  inductively define $s_k= \max\{n\leq N: \om_{n}= \hat\om_{k}\}$, $\hat \om _{k+1}= \om_{s(k)+1}$. If $\om$ is a path in $G$, namely $(\om_{k-1}, \om_{k}) \in E$ for $k=1,\ldots, N$, then
LE$(\om)$ is a self-avoiding path in $G$. For our present purpose we consider the loop erasure LE$(\om^-)$, where $\om^- =(\om_N,\ldots,\om_0)$, the time-reversal of $\om$:\, LE$(\om^-)$ is  obtained from $\om$ by erasing  the loops  in anti-chronological order and tracing the resulting path in that order.

Let $D \in \cD$, a simply connected grid domain (see Section 4),    $v_0, v_e \in  V(\partial D)$ be two distinct boundary vertesies  and $\psi : D\mapsto \H$ be as in the preceding subsection, so that $\psi(v_0)=\infty, \psi(v_e)=0$. Define $\hat o$ by  $\hat o=\psi^{-1}(i)$ and  put $\rho_D=\mbox{in-rad}_{\hat o}(D)$ as before.
 Let  $S ^x $ be  a natural random walk on $G$ started at $x$ (see Section 2 for detailed description) and suppose that  $H_D(v_0,v_e)>0$ (see Section 4.4 for $H_D$). 
Let $\Ga ^{v_0,v_e}$ denote  an excursion,
a natural random walk path,  in $D$ started at $v_0$
and   conditioned to hit $\partial D$ at  $v_e$, where it is stopped.

 We identify  a path in $G$ with the  curve  obtained by linearly interpolating it, and accordingly  use the same expressions $\Gamma^{v_0,v_e}$,  LE$(\Gamma^{v_0,v_e})$ or the like to denote the   (polygonal)  curves  corresponding  to the random walk  path in $G$  the expressions  originally designate, which abuse of notation will not give rise to any confusion.  It is recalled that a chordal SLE$_2$ curve  in $\H$  is the random curve that generates the Loewner chain in $\H$ whose driving function is $\sqrt{2}\, W(t)$, in particular the curve starts at $0$ a.s. Note that  LE$((\Gamma^{v_0,v_e})^-)$ starts at $v_e$, which $\psi$ sends to $0 \in \partial \H$ .
  
%Thm5.5} \label{main1}
\begin{Thm} \label{main1} 
If  the grid domain $D\in \D$ expands to the whole complex plane $\C$ in such a way that $\rho_D \to \infty$ and $D/\rho_D$ is confined in a compact set, then
the simple curve  $\psi \circ {\rm LE}((\Gamma^{v_0,v_e})^-)$  converges to the chordal SLE$_2$ curve  in $\H$
with respect to  driving function.
 \end{Thm}

  The convergence  \lq\lq with respect to  driving function'' in Theorem \ref{main1}  is paraphrased in more precise terms  as follows:   
the driving function of  the random curve  $\psi \circ {\rm LE}((\Gamma^{v_0,v_e})^-)$  under  the conditional probability given that 
 the walk $S^{v_0}$ exits $ D$ through $v_e$  converges weakly (in the usual sense) to     $(\sqrt{2}\, W^{0}(t):t\geq 0)$,  the linear Brownian motion path started at $0$ and  scaled by $\sqrt 2$.   It is noted that the law of the inverse image by $\psi$ of SLE$_2$ in $\H$ is independent of the choice of $\psi$ apart from the time-change  by scaling with a constant factor.

Suzuki \cite{Su} proves the theorem above under the additional assumption that $\partial D$ is locally analytic at $v_0$
and  $S^x$ satisfy invariance principle uniformly for the initial point $x\in V$. As is remarked in \cite{Su} the first assumption 
of the local analiticity at $v_0$  is  imposed to assure the conclusion of Proposition \ref{prop2}
, while the  uniformity of invariance principle with respect to  initial points of the random walk  may be replaced by   our hypothesis  (H) owing to the results of \cite{YY} that are cited as Propositions \ref{YY4.1} and \ref{prop5.3} in the present  paper.  

\v2
  In what follows  we present main steps of the proof of  Theorem \ref{main1}, which is similar to those found in  \cite{SS} or \cite{Su},  focusing  the  description on the new  ingredients that are special to the present setting.
 
We designate  the self-avoiding path LE$((\Gamma^{v_0,v_e})^-)$ in $G$ by
$\gamma =(\gamma _0, \gamma _1, \dots , \gamma _\ell)$ where $\ga_0=v_e$ and $\ga_\ell=v_0$.  By the same symbols    $\ga$ and $\gamma [0,j]$  we also denote  the corresponding  simple curves according to the convention  mentioned previously.

Let $g_t$ and $U(t)$ be  the Loewner chain generated by the simple curve $\psi (\gamma )\subset \H\cup\{0\}$ and its driving function, respectively. 
On recalling  that $g_t$ is defined for the curve $\tilde \ga =\psi(\ga)$ whose time parameter is changed by a function  $\chi$ so that    $2t={\rm hcap}\, \tilde\psi(\ga([0,\chi(t)])$, it is appropriate to   define for $j=0, 1, 2,\ldots$,  $t_j = \chi^{-1}(j)$, namely
$$t_j=\frac{1}{2} {\rm hcap} \,\psi (\gamma [0,j] ).$$
We put
$U_j:=U(t_j)$ and  $D_j:=D \setminus \gamma [0,j]$ and often write $t(j)$ for $t_j$. 
  Following \cite{SS}
we bring in the moving reference point $\hat o_j :=\psi _{t(j)} ^{-1}(i+U_j).$ (Here $\psi_t$ is defined by (\ref{psi_t}), and   $\hat o_j=\hat o(t_j)$ in the notation of Section 5.2.)
In radial case, such a  point  is fixed at the \lq origin'  of $D$ toward which  the loop erasure  evolves, while  in chordal case, $\hat o_j$ must be appropriately  moved along with  $j$  so that there remains a sufficient space around $\hat o_j$ in each  $D_j$, $j=0,1,2\ldots, k$ (for a suitable $k << \ell $), which allows us to apply the invariance principle and its consequences obtained in Section 4.
%which constitutes
%a  sequence of reducing  domains formed by encroachment of  subarcs of $\ga$ into $D$. 

For any $\epsilon >0$,
let 
\[
m:=\min \{ j \geq 1 : t_j \geq \epsilon ^2  \text{ or }  |U_j-U_0| \geq \epsilon \}.
\]

%%%%Lem5.6
\begin{Lem}\label{fund_lem}
There exists a constant $c>0$  such that for each  $\epsilon >0$ and $M>1$,
there exists $R>0$ such that
if $D\in \mathcal D_{R,M}$, then
\[|U_m-U_0| <2\e,\quad 
|\E [ U_m-U_0]| \leq c \epsilon ^3
\quad\mbox{and} \quad
|\E [(U_m-U_0)^2-2t_m ]| \leq c \epsilon ^3.
\]
(Although $U_0=0$, we enter $U_0$ in the formulae above to indicate how they show   when the starting position 
 $\gamma _0$ of the curve $\ga$ is not mapped to the origin by $\psi$.)
\end{Lem}
%\begin{proof}

Proving this lemma is an essential step for  Theorem \ref{main1}.  We give a proof, omitting some details,   and  indicate where we need   Proposition \ref{prop2}.  
We make use of  a martingale which is suggested  in \cite{LSW} and adopted in \cite{YY}, \cite{Su}  as an observable 
for the same purpose as ours.    We define
\beqn\label{M_j}
M_j:=\frac1{Z_0} \cdot\frac{H_{D_j}(w, \gamma _j)}{H_{D_j}(v_0, \gamma _j)},
\eeqn
(for any $\delta >0$ and $w \in V(D)$), where 
$$Z_0= Z(D, \hat o,v_e) =  \frac{H_{D}(\hat o, v_e)}{H_{D}(v_0, v_e)}.$$
Then $M_j$ is a martingale with respect to the filtration generated by $ \ga[0,j]$, $j=0,1,2,\ldots$, of which fact an abridged proof is given   in \cite{LSW} (see  \cite{YY}, \cite{Su} for a detailed proof). 

\v2
{\sc Remark 4.} $Z_0$ is unbounded and acts  as a  normalizing  constant.  In \cite{LSW} and \cite{YY},  the radial SLE  being  concerned,  $v_0$ is replaced by $\hat o$ and  the normalization is not needed (one may take $Z_0=1$).  Suzuki \cite{Su}, dealing with the chordal case, adopts a different normalization: $Z_0 = H_D(v_0; A)$, where $A=\psi^{-1}([-1,1])$.  The difference is of only technical matter, although with  our choice the proof is simpler (owing to Corollary\ref{prop5.2}).   

\v2
{\it Proof   of Lemma \ref{fund_lem}.}   \; First we notice a  fact  which  underlies  the arguments  given below.  It is shown  in \cite{SS} (the first paragraph of the proof of Proposition 4.1 of it)  that there exists a universal constants $c$ and $R$ such that if $\rho_D\geq R$ and $\e<1/R$,
 \beqn\label{ser}
 \mbox{in-rad}_{\hat o}(D_m) \geq c \rho_D.
 \eeqn 
(In \cite{SS} the boundedness of edge length over $G$ is used, which can be plainly  replaced by the property (1) of Section  4.1.)

  We are to derive a neat expression of  $1/Z_0 M_m $ as well as $1/Z_0$ (see  (\ref{A2}) and (\ref{A1})  below)  by applying  Proposition \ref{prop5.3} (Poisson kernel approximation).
Let  $0<\epsilon <1/2 $  and  let $B_{r}(a) $ and $U_\de$ be defined  as in Section 3.2 (with  $(D, \hat o)$). Then, on writing $B(\e^2, v_0) $ for  $B_{\e^2}(v_0) $ 
\beqn\label{H_dcmp}
H_{D_j}(v_0, \gamma _j)
=\sum _{y \in C_\e} P [ S^{v_0} (\tau _{B(\e^2,v_0)})=y] H_{D_j}(y, \ga_j). %, S^{v_0} ( \tau _{D_j}) = \gamma _j),
\eeqn
where $C_\e$ is the set of all vertices $y$ for which the summand is positive (hence located along $D \cap \partial  B_{\e^2}(v_0))$. 
We split the sum on the right-hand side  into two  parts according as  $y\in U_\de$ or not. It is here  that we apply  Proposition \ref{prop2}. Owing to it we can choose $\de=\de(\e)>0$ (independent of $D_j$) so that for  $D$ with $\rho_D$ large enough, 
 \[
P \Big[ |\fa_D(S^{v_0} (\tau _{B(\e^2, v_0)}) )| >1-\de \,\Big|\, S^{v_0} ( \tau _{B(\e^2, v_0)}) \in D\Big] =O(\epsilon ^{3}).
\]
It is easy to see that $\ga[0,t_m]$ is confined in a small neighborhood of  $v_0$  (cf. the argument given for  (\ref{u7}) below if necessary) and  on using  Corollary \ref{prop5.2} we observe that the conditioning on  the event  $S^{v_0} ( \tau _{B(\e^2,v_0)}) \in D$ can be replaced by conditioning on
$ S^{v_0} ( \tau _{D_j}) = \gamma _j$, 
 implying that on the right-hand side of   (\ref{H_dcmp}) the  proportion of  the contribution of the sum on $y\in U_\de$ to the whole sum is $O(\epsilon ^{3})$.

Once $\de$ is determined  we apply Proposition \ref{prop5.3} or rather Corollary \ref{cor5.30}  for computation of  the sum over $V(D)\setminus U_\de$. It is immediate to see
\beqn\label{A2}
\frac1{Z_0}
= \sum_{y\in  C_\e \setminus U_\de} \frac{\Im\, \psi  (y)}{|\psi  (y) - U_0|^2} \Big(1+O(\epsilon ^3 )\Big)
\eeqn
(provided that $\rho_D$ is large enough depending only on $\e$).
We wish to apply  Corollary \ref{cor5.30}  with  $D_m$ in place of $D$, which is to result in the formula
\beqn\label{A1}
\frac1{Z_0M_m} =\frac{H_{D_m}(v_0, \gamma _m)}{H_{D_m}(w, \gamma _m)} 
= \frac{|\psi_{t(m)} (w) -U_m|^2}{\Im\, \psi_{t(m)} (w)} \sum_{y\in C_\e \setminus U_\de} \frac{\Im\, \psi_{t(m)}(y)}{|\psi_{t(m)}(y)-U_m|^2} 
 \Big(1+O(\epsilon ^3 )\Big).
\eeqn
 Since $\psi(\partial B_{\e^2}(v_0) \setminus U_\de)$  is distorted by the mapping $g_{t_m}$, this application requires justification (specifically to ensure the condition that   $\psi_{t_m}(y)$ and  $\psi_{t_m}(w)$   remain in a compact set  of $\H$), which however is given by Corollary \ref{cor_k-est}  (see (\ref{K}) below).

In view of (3) of Section 4.1 we have 
$${\rm diam} (\psi_{t(m-1)} (\gamma [t_{m-1},t_m]) ) \leq \epsilon^2 $$
whenever $\rho_D$ is large enough, which fact together with  Lemma \ref{k-est2} (its harder half) applied to the Loewner chain  $\tilde g_t :=g_{t+t(m-1)}\circ g_{t(m-1)}^{-1}$  implies $t_m-t_{m-1} =O(\e^4)$ and $U(t_m)- U(t_{m-1}) =O(\e^2)$; hence, by the definition of $m$,
\begin{equation}
 t_m <2\e^2\quad \mbox{and} \quad \sup\{|U(s)-U(0)|: s\in [0, t_m]\}  <2\e.    \label{u7}
\end{equation}

%%%
For  $z\in \D$ we have  $$\psi \circ \fa^{-1}_D(z)=i\frac{1-z}{1+z},$$
in particular   for $y\in C_\e$,  $|\psi(y)| \sim 2/\e^2$.
%%%%%
Hence,  according to  Corollary \ref{cor_k-est}  we obtain for $y\in C_\e$
\begin{equation}\label{K}
|\psi_{t(m)} (y)-\psi (y)|  = O(\epsilon ^4) \quad \text{and} \quad \frac{\Im\, \psi_{t(m)}(y)}{\Im\, \psi(y)} =  1 +  O(\e^6);
\end{equation}
moreover,  using   the Loewner equation (\ref{cle}) in addition, we also  infer that
if $w \in V(D)$  is subject to 
 \begin{equation} \label{u00}
\Im\, \psi (w) \geq \frac{1}{2} , \  |\psi (w)| \leq 3,
\end{equation}
then
\beqn\label{Leq}
\psi_{s} (w) -\psi (w) = \frac{2s}{\psi(w)-U_j} + O(\epsilon ^3) \quad \mbox{   for}  \quad s \in [0,t_m],
\eeqn %\label{est_w}
%\end{equation}
entailing
$
\Im\, \psi_{t(m)} (w) \geq \frac{1}{3},\  |\psi_{t(m)}(w)|\leq 4$.
 % $\frac12 \la \leq |\psi_{t(m)}(y)| \leq  3\lambda$.
Noting     $U_m/\psi(y) =O(\e^3)$ we apply  the relation  (\ref{K}) to find   that the ratio of each term of the sum in (\ref{A2}) and the corresponding one of the sum appearing  in (\ref{A1})  is $1+O(\e^3)$, hence so is  the  ratio of  the two sums.
Hence by  (\ref{A1}) and   (\ref{A2})
\begin{align*}
M_m=  \,&\frac{\Im\, \psi_{t(m)}(w)}{|\psi_{t(m)}(w)-U_m|^2} \Big(1 +O(\epsilon ^3 )\Big)\\
= \, & \Im\, \left (\frac{-1}{ \psi_{t(m)}(w)-U_m}\right )+O(\epsilon ^3). 
\end{align*}

Let $F(z,u)=\Im [-1/ (z-u)], z\in\C, u\in \R$ and compute the difference $F(\psi_{t(m)}(w),U_m)- 
F(\psi(w), U_0)$ by means of  Taylor expansion. With the help of  (\ref{Leq}) as well as $0=E[M_m-M_0] =  E[F(\psi_{t(m)}(w), U_m) - F(\psi(w),U_0)]$, this leads to
\begin{equation} \label{u9}
\Im\left \{ \frac {1}{(\psi_{t(m)} (w)-U_0)^2} E [U_m-U_0]+
\frac {1}{(\psi_{t(m)} (w)-U_0)^3} E [(U_m-U_0)^2-2t_m]\right\} 
 =O(\epsilon ^3). \quad
\end{equation}
 Take two distinct vertices $w_1$ and $w_2$ in place of $w$ so that $\psi(w_1) =i +O(\e^3)$ and  $\psi(w_2)= e^{i\pi/3} +O(\e^3)$, and you find out the required  relation of the lemma. \qed

%%%%%
%%%%
\v2
{\it  Proof of Theorem \ref{main1}}.  
Having proved  Lemma \ref{fund_lem}  it is easy to adapt the arguments given in \cite{LSW} for our proof of  Theorem \ref{main1}.
 Let $T>1$ and $\epsilon _1 >0$  and put  $\tilde T = \sup\{t\in [0,T]:|U(t)|<1/\epsilon _1\}$. 
Let $\epsilon >0$ be small enough.
Let $m_0=0$ and define $m_n$ inductively by
\[
m_n := \min \{j > m_{n-1} : t_j- t_{m_{n-1}} \geq \epsilon^2 \text{ or } |U_j - U_{m_{n-1}}|\geq \epsilon \}.
\]
Let $N:= \max \{ n \in \N : t_{m_n} < \tilde{T} \}$.
  By a   Markovian nature of  the walk $S^{v_0}$ conditioned on $\ga[0,j]$ (Lemma 3.2 of \cite{LSW})   together with the Huygens property (\ref{U2}) of the Loewner chain, we  apply Lemma  \ref{fund_lem} with $(D_{m_n}, v_0,\gamma _{m_n}, \hat o_{m_n}, U(\cdot+t_{m_n})$ in place of $(D,v_0, v_e, \hat o, U)$. The application is secured   owing to Lemma \ref{refpt}, which shows that $\mbox{in-rad}_{\hat o_j}(D_j)/\rho_D$ is bounded below by a positive constant for $n\leq N$. It then follows that
  there exists $R=R(\epsilon, \epsilon _1, T, M )>0$ such that if $D\in \cD_{R,M}$, then for any $n \leq N$
\[  | U_{m_{n+1}}-U_{m_n}| <2\e,\quad
\E [ U_{m_{n+1}}-U_{m_n}\, |\, \gamma [0, m_n] \,] =O(\epsilon ^3),
\]
and
\[
\E  [(U_{m_{n+1}} -  U_{m_n})^2\ | \ \gamma [0, m_n] \,] 
=\E [2( t_{m_{n+1}} -t_{m_n}) \,|\,  \gamma [0, m_n] \,] +O(\epsilon ^3).
\]
 It is a more or less  standard issue   of the probability theory  to deduce from these three relations  that the law of $U(t)$ weakly converges to that of
 a scaled Brownian motion $\sqrt{2}\, W(t)$ as $\rho_D\to\infty$ and is carried out in \cite{LSW} (Section 3.3: especially the arguments following  Eqs.(3.16-17)).  \qed
 
  \v2\v2

%to  find that 
%$\gamma ^{(t_{m_n})}(\cdot )=\gamma ( t_{m_n}+ \cdot )$ is the same distribution as 
%the time-reversal of the loop erasure of 
%a natural random walk on $G_\delta$ started at $b_\delta$ and stopped on exiting $D_{m_n}$ 
%and conditioned to hit $\partial D_{m_n}$ at $\gamma _{m_n}$.
%We apply Lemma  \ref{fund_lem} with $(D_{m_n}, \gamma _{m_n}, b, p_{m_n})$ for any $n \leq N$.
%Then, we deduce from the fact stated at (\ref{U2}) that 

{\sc Remark 5.} {\it Another observable.}\;  In \cite{LSW} the random variables 
$h^+_j = $the number of visits to $w$ by $\Ga^{v_0,\ga_j}[0,n_j]$ is adopted as a martingale observable  for deriving the estimate corresponding to Lemma \ref{fund_lem} (see Proposition 3.4 of \cite{LSW}). In the present setting we have
  $$E\Big[h_j^+\, \Big|\,\ga[0,j]\Big]=\frac{\hat G_j(v_0,w)\hat H_j(w,\ga_j)}{\hat H_j(v_0,\ga_j)},$$
  where $\hat\,$ indicates the corresponding objects for the walk conditioned to exits $D$ through $v_e$.
  We can adapt the proof in \cite{LSW}, provided that in our assumption (H) the weak convergence is
 understood  relative to not  the metric $d_{\D}$ but the usual  metric for uniform convergence of functions of $t\geq0$  on each finite interval.  This modification of (H) is needed since for the observable we need to approximate the Green function $G_j(y,w)$ by the corresponding one for Brownian motion.
\v2
{\sc Remark 6.} \, From the convergence of the driving process it follows that the image of loop erasure $\psi$(LE$((\Ga^{v_0,v_e})^-))$ converges relative to Hausdorff metric if the time  (parametrized by the h-capacity)  is restricted to a finite interval (cf. \cite{LSW}, \cite{SS}).

\v2\v2\n
{\sc 5.4. \; Uniform convergence.}
\v2
 Suzuki \cite{Su} points out  that 
   if  the transition probability $p(u,v)$ has an invariant measure $\pi$  and the dual walk with respect to  it as well as the original walk 
   satisfies invariance principle (as given in \cite{Su}),  then  the convergence of $\psi$(LE$(\Ga^{v_0,v_e}))$ with respect to the metric of  \lq uniform'  convergence of path functions    follows from their convergence  with respect to  driving function (i.e., convergence of their  driving functions). His reasoning is  based on a result by Shefield and Sun \cite[Corollary 1.7]{sSS}, which asserts that 
 if a sequence of random simple curves  in $\D$ and that of their time-reversals both converge to SLE$_{\k}$ curves with $\k\leq  4$ with respect to  driving function,  then the weak convergence with respect to the metric $d^*_\D$ holds true.  
  It is shown by Lawler  \cite{L2} that for symmetric random walks (in fact Markov chains),  the law of LE$((\Gamma^{v_0,v_e})^-)$ agrees with that of  (LE$(\Gamma ^{v_0,v_e}))^-$.  Suzuki observes   that  the symmetry assumption is dispensable so that  the result  extends   to general Markov chains. 
   Note that on using  this  result,  [LE$((\Ga^{v_0,v_e})^-)]^-$ has the same distribution as $L(([\Ga^*]^{v_e,v_0})^-)$ under the existence of $\pi$, where $\Ga^*$ designates the excursion of the dual walk. (Cf. Section 5.2 of \cite{Su} for more details). In below we formulate the corresponding result in the usual setting in which we consider the loop erasure of the scaled  random walk    that is confined  in a fixed domain. For simplicity we shall assume
    $p(u,v)= p(v,u)$ so that the  condition on the dual walk mentioned above is plainly true.  
 
 %{\sc Boundary correspondence via a conformal map.} ~ 
 Let $D$ be a bounded and simply connected domain of $\C$ (not necessarily a grid domain). We  are concerned with the correspondence between points of $\partial D$ and those of  the unit circle  $\partial \D$, which is nicely given    by means   of {\it prime ends} of $D$,  ideal boundary points.  We do not give any definition of prime ends   for which the readers may be referred to \cite{Pm} or \cite{A}.  %Instead, we state  a fundamental fact concerning it.  
 What is  needed in this paper are the facts mentioned below.    There is a natural one-one correspondence between the set of prime ends and $\partial U$. We write $\partial_{pnd} D$ for the set of all prime ends of $D$ and denote the correspondence by $\iota$.
   Let $\fa: D\mapsto \D$ be any  conformal map and denote by $\fa \cup \iota$  its extension to $D \cup \partial_{pnd}  D$ by means of $\iota$.
  A  topology of  $D \cup \partial_{pnd}  D$ is given so that  $\fa \cup \iota$ is   a homeomorphism   from   $D \cup \partial_{pnd}  D$ onto $\overline \D$. 
For $\zeta\in \partial \D$ 
  let $C(\fa, \zeta)$ be the set of all limit points of $\fa^{-1}{z}$ as $z\in \D \to \zeta$:
 $$C(\fa, \zeta) = \bigcap_{r>0}\overline{\{\fa^{-1}(z): z\in \D, |z-\zeta|<r\}}. $$ 
 The  correspondence $\bar \fa$ is natural in the sense that the topology induced is independent of the choice  of $\fa$ and  if a sequence  $x_n\in D$ converges to $x\in \partial_{pnd}  D$, then  dist$(x_n, C(\fa, \iota(x)))$ tends to zero (and under  an additional condition on $(x_n)$ the converse holds).  (Cf. Section 2.5 of \cite{Pm}.)

%Then each prime end is identified with $C(\fa, \zeta) $  for the  $\zeta$  that corresponds to the prime end  (\cite{Pm}, Theorem 2.16).  

%It follows  from  the very definition that $C(\fa, \zeta)$  is a connected closed subset of $\partial D$ and constitutes a partition of  $\partial D$.  

\v2
% {\sc  Convergence of LERW to SLE$_2$.} ~
%Let $D$  be a simply connected domain,  $\fa$ be a conformal map of $D$ onto $\D$ as above and  $a, b \in %\partial_{pnd} D$  be two distinct prime ends.
%Let $G=(V,E)$ be a planar graph as described above and let $S^v= (S_n^v)_{n=0}^\infty$ be a random walk  defined on some probability space $(\Om, \F, P)$  started at $v\in V$ with a transition probability $p(u,v)$ (see Section 3 for detailed description).  For $w\in V$, let  $\sigma_w$  be the first passage time of  $w$ by the random walk that is under consideration.
%Let $\Ga ^{u,v}$ denote  an excursion, a natural random walk path,  in $D$ started at $u$ and at  $v$, where it is stopped.
%If $D$ is a Jordon domain,  the following theorem asserts that   $S^{v_0}$  converges to a chordal SLE$_2$ in a 
%reasonable sense. 

\begin{Thm} Suppose, in addition to the basic hypothesis $(H)$,  that $p$ is symmetric: $p(u,v) = p(v,u)$.  Let $D$  be a simply connected domain,  $\fa$ be a conformal map of $D$ onto $\D$ (as above) and  $a, b \in \partial_{pnd} D$  be two distinct prime ends of $D$.
 Let  $(\rho_n)$ be  a sequence of positive numbers and $(u(n))$ and $(v(n))$  be two sequences in $V(D)$    such that    both $u(n)/\rho_n$ and  $v(n)/\rho_n$ are in $\bar D$ and  $\rho_n \to\infty$,    $u(n)/\rho_n  \to a$ and $v(n)/\rho_n\to b$ as $n\to \infty$ and that   $S^{u(n)}$,  the random walk started at  $u(n)$,  arrives at $v(n)$ before exiting $\rho_n D =\{\rho_n x: x\in D\}$ with positive probability. Let $\mu_n$ be the conditional probability law of the path $\Ga_n := (S^{u(n)}(k): k=0,\ldots, \sigma_{v(n)})$, %the path of the random walk  $S^{a(n)}_k$ for $k=0,\ldots, \sigma_{b_n}$ 
 given that $S^{u(n)}$  visits   $v(n)$ before exiting $\rho_nD$,  and $\ga_n$ be  the linear interpolation of the loop erasure $LE(\Ga_n)$.  Then  the  law of $\fa\circ \ga_n$  induced from $\mu_n$ weakly converges relative to the metric $d_{\D}$  to that of  the  chordal SLE$_2$ curve from $\fa(a)$ to $\fa(b)$ in $\D$.
\end{Thm}
\v2\n
\pf\;  We reduce the problem   to Theorem \ref{main1}.  To this end let $V_n$ be the set of all vertices  that  the walk $S^{(u(n)}$  arrives before exits  $\rho_n D$ with positive probability:
$$V_n = \{ w\in V: P[ S^{u(n)}_k =w\,\,   \mbox{for some}  \,\, k< \tau_{\rho_n D}]>0\}$$
and  define a grid domain  $D_n\in \cD$ as  the smallest one among those that  contains $V_n$ (in their interior). Then 
   $S^{u(n)}$ exits  $D_n$ and  $\rho_n D$ at the same time and  $D_n/\rho_n$ converges to $D$ in the Carath\'eodory sense (see   \cite[Proposition 3.63]{L}).   Noting that  in Theorem \ref{main1}  $v_0$ and/or $v_e$ have not to be taken  from  boundary; they may be vertices near the boundary (see Proposition \ref{prop5.3}), the premise  in    Corollary 1.7   of \cite{sSS}  mentioned above is readily verified. \qed

   \v2\v2\n
 % {\bf Acknowledgement.}  I  express my gratitude to Izumi Okada, who kindly produced the figures I through III. 
   % I would like to thank Ryoki Fukushima for   pointing  out that  one needs Lemma \ref{refpt} in order to follow the arguments of \cite{SS} in which no reference  to  the corresponding assertion is made. 
%%%%%%%%
%%%%%%%%
%and it follows from Lemma \ref{refpt} that ${\rm rad}_{p_m}(D_m) \geq r^\prime$   for some $r^\prime>0$.
%Therefore, we can apply Proposition \ref{est_mart} to the domain $D_m$ with the reference point $p_m$,

\end{document}